\documentclass[12pt]{article}
\usepackage{color}

\usepackage{float}
\usepackage{tikz} 
\usepackage{amsmath}

\usepackage{hyperref}
\hypersetup{hypertex=true,
	colorlinks=true,
	linkcolor=blue,
	citecolor=blue}

\usetikzlibrary{shapes,snakes}

\newcommand{\qed}{\hfill {$\diamond$}}
\newcommand{\pf}{{\bf Proof: }}

\newcommand{\2}{\vspace{0.2cm}}
\newcommand{\dom}{\rightarrow}

\newcommand{\compdom}{\Rightarrow}
\newcommand{\ncompdom}{\not\Rightarrow}
\newcommand{\smd}{semicomplete multipartite digraph}

\newtheorem{theorem}{Theorem}[section]
\newtheorem{proposition}[theorem]{Proposition}
\newtheorem{conjecture}[theorem]{Conjecture}
\newtheorem{lemma}[theorem]{Lemma}
\newtheorem{claim}{Claim}
\newtheorem{corollary}[theorem]{Corollary}

\newtheorem{problem}[theorem]{Problem}
\newcommand{\yw}[1]{{\color{red}#1}}

\begin{document}
	\title{Generalized paths and cycles in semicomplete multipartite digraphs}
	\author{J\o rgen Bang-Jensen\thanks{
			Dept.\ of Math.\ and Compt.\ Sci., University of Southern Denmark,
			Odense, Denmark (email: jbj@imada.sdu.dk).}\and
		Yun Wang\thanks{
			School of Mathematics, Shandong University, Jinan 250100, China and Dept.\ of Math.\ and Compt.\ Sci., University of Southern Denmark, Odense, Denmark (email: yunwang@imada.sdu.dk,wangyun\_sdu@163.com).}\and
		Anders Yeo\thanks{
			Dept.\ of Math.\ and Compt.\ Sci., University of Southern Denmark, Odense, Denmark 
			and Department of Mathematics and Applied Mathematics, University of Johannesburg, Auckland Park, 2006 South Africa (email: yeo@imada.sdu.dk).}}
	\maketitle 
	
	\begin{abstract}
		A digraph is {\bf semicomplete} if it has no pair of non-adjacent vertices. It is {\bf complete} if every pair of distinct vertices induce a 2-cycle.
		It is well-known and easy to show that even the following version of the directed travelling salesman problem is NP-complete: Given a strongly connected complete digraph $D=(V,A)$, a cost function $w: A\dom \{0,1\}$ and a natural number $K$; decide whether $D$ has a directed Hamiltonian cycle of cost at most $K$. We study the following variant of this problem for $\{0,1\}$-weighted semicomplete digraphs where the set of arcs which have cost 1 form a collection of vertex-disjoint complete digraphs. A digraph is \textbf{semicomplete multipartite} if it can be obtained from a semicomplete digraph $D$ by choosing a collection of  vertex-disjoint subsets $X_1,\ldots{},X_c$ of $V(D)$ and then deleting all arcs both of whose end-vertices lie inside some $X_i$. Let $D$ be a semicomplete digraph with a cost function $w$ as above, where $w(a)=1$ precisely when $a$ is an arc inside one of the subsets
		$X_1,\ldots{},X_c$ and let $D^*$ be the corresponding \smd{} that we obtain by deleting all arcs inside the
                $X_i$'s. Then every  cycle $C$ of $D$ corresponds to a {\bf generalized cycle} $C^g$ of $D^*$ which is either the cycle $C$ itself if $w(C)=0$ or a collection of two or more paths that we obtain by deleting all arcs of cost 1 on $C$. Similarly we can define a {\bf generalized path} $P^g$ in a semicomplete multipartite digraph. The purpose of  this paper is to  study structural and algorithmic properties of generalized paths and cycles in semicomplete multipartite digraphs. This allows us to identify classes of directed $\{0,1\}$-weighted TSP instances that can be solved in polynomial time as well as others for which we can get very close to the optimum in polynomial time. Along with these results we also show that two natural questions about properties of cycles meeting all partite sets in semicomplete multipartite digraphs are NP-complete.\\

                \noindent{}{\bf Keywords:} Semicomplete multipartite digraph; Generalized cycle; $\{0,1\}$-TSP problem; Hamiltonian cycles; NP-completeness
		   
	\end{abstract}
	
	\section{Introduction}
	A digraph $D=(V,A)$ is {\bf semicomplete multipartite} if its vertex set has a partition $V=V_1\cup \cdots{}\cup{}V_c$ so that $A$ contains no arc with both ends in some $V_i$ and $A$ contains at least one arc between $u,v$ for every choice of $i,j\in [c]$ where $i\neq j$ and $u\in V_i,v\in V_j$. The sets $V_1,\ldots{},V_c$ are called the {\bf partite sets} of $D$. A special well-studied subclass of semicomplete mulitpartite digraphs is the class of  {\bf extended semicomplete digraphs} which are those   semicomplete multipartite digraphs $D=(V_1,\ldots{},V_c,A)$ for which we have $V_i\compdom V_j, V_j\compdom V_i$ or $D[V_i\cup V_j]$ is a complete bipartite digraph for every choice of $i,j\in [c]$ such that $i\neq j$.

	Let $D=(V_1,\ldots{},V_c,A)$ be a semicomplete multipartite digraph. A {\bf generalized cycle} 
	 of $D$ is a collection of vertex-disjoint directed paths $P_{u_1,v_1},P_{u_2,v_2},\ldots{},P_{u_r,v_r}$ of $D$ with the property that $v_i$ and $u_{i+1}$ belong to the same partite set of $D$ for $i\in [r]$, where we take $u_{r+1}=u_1$. A {\bf generalized path} is defined similarly except that now we do not require that $v_r$ and $u_1$ belong to  the same partite set of $D$ (but they may do so, in which case they must be distinct). To be more precise we define  the generalized path $P_{u_1,v_1},P_{u_2,v_2},\ldots{},P_{u_r,v_r}$ to start in $u_1$ and end in $v_r$
	and  call it a generalized $(u_1,v_r)$-path. For simplicity we call a generalized cycle, respectively a generalized path a {\bf G-cycle}
	respectively a {\bf G-path}. The {\bf length} of a G-cycle $C^g$ or G-path $P^g$, denoted by $\ell{}(C^g)$ and $\ell{}(P^g)$, respectively,
	is the number of arcs of $A$ belonging to the G-cycle or G-path. Note that a G-cycle or G-path may have length 0, in which case all its vertices belong to the same partite set of $D$. It follows from the definition of length that a G-cycle of length $n=|V(D)|$ is a Hamiltonian (real) cycle of $D$.
	
	A G-path or G-cycle in a \smd{} $D=(V_1,\ldots{},V_c,A)$ is  {\bf good} if it contains at least one vertex from each partite set $V_i$, $i\in [c]$ and it is {\bf spanning} if it contains all vertices of $D$. Note that when $D$ is a semicomplete digraph, then the only G-cycles and G-paths $D$ can have are directed cycles and paths in $D$. Hence, as already indicated in the abstract, (good) G-cycles and G-paths generalize (Hamiltonian) paths and cycles in semicomplete digraphs. Furthermore, by the remarks in the abstract, the problems of finding a Hamiltonian path or cycle of minimum cost in a semicomplete digraph $D$ with a $\{0,1\}$-cost function $c$ on the arcs so that the arcs of cost 1 form vertex-disjoint complete subdigraphs of $D$ are equivalent to the following problems on G-paths and G-cycles, respectively.

	\begin{problem}\label{prob:longestG-path}
		Given a semicomplete multipartite digraph $D=(V,A)$. Find a spanning G-path which uses the maximum number of arcs of $A$.
	\end{problem}
	
	\begin{problem}\label{prob:longestG-cycle}
		Given a semicomplete multipartite digraph $D=(V,A)$. Find a spanning G-cycle which uses the maximum number of arcs of $A$.
	\end{problem}

        Note that while the longest path problem in semicomplete multipartite digraphs has been completely solved \cite{gutinK1988,gutinSJDM6}, the complexity of the longest cycle problem in \smd{}s is still open.

        The purpose of this paper is to  study  the two problems above as well as related questions regarding G-paths and G-cycles in semicomplete multipartite digraphs. Among other results we  show that Problem \ref{prob:longestG-path} can be solved in polynomial time and obtain a number of results on Problem \ref{prob:longestG-cycle}. We also indicate why Problem \ref{prob:longestG-cycle}, which is still open for the full class of semicomplete multipartite digraphs,  is much harder than Problem \ref{prob:longestG-path}.\\

          The paper is organized as follows. In Section 2 we give necessary terminology, some basic lemmas on G-cycles and G-paths as well as a number of results on semicomplete (multipartite) digraphs that will be used in our proofs. In particuilar we show how to obtain a G-cycle factor with the maximum number of arcs in a \smd{} in polynomial time. In Section 3 we prove some
          basic results about G-paths and G-cycles which illustrate the diffence between these objects and 'real' paths and cycles. In Section 4 we use results on hamiltonian paths with prescribed end vertices in semicomplete digraphs to derive analogous results for G-paths in semicomplete multipartite digraphs. In Section 5 we solve Problem \ref{prob:longestG-path} and show that the characterization that we obtain is closely related to the characterization in \cite{gutinK1988,gutinSJDM6} for longest paths in semicomplete multipartite digraphs. In Section 6 we solve Problem \ref{prob:longestG-cycle} for extended semicomplete digraphs, a well studied subclass of semicomplete multipartite digraphs.\\
          The key to obtaining the polynomial algorithm for the hamiltonian cycle problem for \smd{}s
          in \cite{bangJGT29} was the study of so-called irreducible cycle factors in \cite{yeoJGT24}. In Section 7 we study the analogous problem of finding an irreducible  G-cycle factor containing  the maximum number of arcs in \smd{}s. Our proofs imply that there is a polynomial algorithm for converting any G-cycle factor $F$ into an irreducible G-cycles factor $F'$  which has at least as many arcs as $F$. If this  G-cycle factor has only one G-cycle, then that G-cycle is a solution to Problem \ref{prob:longestG-cycle}. In general an irreducible factor may have many G-cycles. In that case, similar to the results in \cite{yeoJGT24}, we obtain a lot of information about the structure of the \smd{} induced by the vertex set of such an irreducible  G-cycle factor. In Section 8 we prove a number of results concerning Problem \ref{prob:longestG-cycle} and illustrate by several examples why the problem is difficult. Recall again that also the longest cycle problem in \smd{}s is open. In Section 9 we translate our results to results for $\{0,1\}$-weighted TSP instances  on complete digraphs. Finally in Section 10 we show that two natural questions about cycles visiting all partite sets in \smd{}s are NP-complete.

	\section{Terminology and Preliminaries}
	
	A digraph  $D=(V,A)$ is {\bf strong} if there exists a path from $x$ to $y$  in $D$ for every ordered pair of distinct vertices $x$, $y$ of $D$ and $D$ is {\bf $k$-strong} if $|V|\geq k+1$ and  $D-X$ is strong for every subset $X\subset V$ of size at most $k - 1$. An {\bf $(x,y)$-path} is a path starting in $x$ and ending in $y$. If there is an arc $xy$ in $D$, then we say that \textbf{$x$ dominates $y$} and use $x\dom y$ to denote it.
	A {\bf $k$-path-cycle subdigraph} of a digraph $D$ is a collection of $k$ paths $P_1,\ldots{},P_k$ and zero or more cycles $C_1,\ldots{},C_r$ of $D$ such that all of these paths and cycles are vertex-disjoint. A 0-path-cycle subdigraph is called a {\bf cycle subdigraph} and if it covers all vertices of $D$, it is a {\bf cycle factor} of $D$. Similarly we can define a {\bf $k$-G-path-G-cycle subdigraph} and a {\bf G-cycle factor} by just replacing path by G-path and cycle by G-cycle above. In particular, a
	{\bf G-cycle factor} of a semicomplete multipartite digraph $D$ is a vertex-disjoint collection of G-cycles
	$C^g_1,\ldots{},C^g_k$, $k\geq 1$ such that every vertex of $D$ is contained in exactly one $C^g_j$, $j\in[k]$.
	Recall that we allow a G-cycle $C^g$ to contain no arcs of $D$, in which case all vertices of $C^g$ belong to the same partite set of $D$.

	Let $L=v_1v_2\cdots v_{l}, l\geq 2$ be a (G-)path or a (G-)cycle.  The successor of $v_i$ on $L$, written as $v_i^+$, is the vertex $v_{i+1}$. Similarly, one can denote the predecessor of $v_i$ by $v_i^-$. We denote $L[v_i,v_j]=v_iv_{i+1}\cdots v_j$. For disjoint subsets $X,Y$ of $V(D)$, a (G-)path $P$ is an $(X,Y)$-(G-)path if it starts at a vertex $x$ in $X$ and ends at a vertex $y$ in $Y$ such that $V(P)\cap(X\cup Y)=\{x,y\}$.\\
	
	The following result follows from the fact that, using e.g. minimum cost flows, we can find a cycle factor of minimum cost in any digraph with non-negative weights on the arcs: Given a semicomplete multipartite digraph $D=(V,A)$ whose partite sets of size at least 2 are $V_1,V_2,\ldots{},V_k$ we let $D^*$ be the arc-weighted semicomplete digraph that we obtain from $D$ by adding all possible arcs between vertices in the same $V_i$, $i\in [k]$ and give these weight 1 and give weight 0 to all arcs in $A$. Now it is easy to see that every minimum cost cycle factor in $D^*$ corresponds (by deleting all arcs inside $V_i$'s) to a G-cycle factor of $D$ that covers the maximum number of arcs of $D$.
	
	\begin{lemma}\label{lem:maxnoarcsCF}
		There is a polynomial algorithm for finding a G-cycle factor containing the maximum number of arcs in a \smd{}. \qed\end{lemma}
	
	Let $D$ be a digraph.  For disjoint subsets $X, Y$ of $V(D)$, we use the notation
		$X\dom Y$ to denote that $x\dom y$ for any pair of adjacent vertices $x\in X$ and $y\in Y$. $X\compdom Y$ means that $X\dom Y$ and no arc from $Y$ to $X$. For two vertex-disjoint subdigraphs $D_1,D_2$ of $D$, we will often use the shorthand notation $D_1\dom D_2$ and $D_1\compdom D_2$ instead of  $V(D_1)\dom V(D_2)$ and $V(D_1)\compdom V(D_2)$.

	We now recall some results that will be used or generalized in this paper.

	Let $P^g$ be an $(x,y)$-G-path in a \smd{} $D$ (possibly $x=y$ and $P^g$ has order one) and let $Q^g=v_1\cdots{}v_t$ be a G-path or a G-cycle in $D-V(P^g)$. Then we say that $P^g$ has a \textbf{partner} on $Q^g$ if there exist consecutive vertices $v_i,v_{i+1}$ of $Q^g$  such that $v_i$ dominates $x$ and $y$ dominates $v_{i+1}$ ($v_iv_{i+1}$ is called the partner of $P^g$). Observe that if $P^g$ has a partner on $Q^g$, then the G-path $P^g$ can be inserted into $Q^g$ to give a new G-path (or G-cycle) $L^g=Q^g[v_1,v_i]P^gQ^g[v_{i+1},v_t]$ such that $\ell{}(L^g)\geq \ell{}(Q^g)+1$.  More precisely we will have $\ell{}(L^g)=\ell{}(P^g)+\ell{}(Q^g[v_1,v_i])+\ell{}(Q^g[v_{i+1},v_t])+2$ and this is at least $\ell{}(P^g)+\ell{}(Q^g)+1$ as we gain the arcs $v_ix$ and $yv_{i+1}$ and loose at most one arc, namely $v_iv_{i+1}$.\\
	
	The following lemma generalizes Corollary 6.4.18 in \cite{bang2009} for semicomplete multipartite digraphs.
	\begin{lemma}\label{partner-spanning}
		Let $D$ be a semicomplete mulitpartite digraph. Suppose that $P^g=u_1\cdots{}u_r$ is a G-path and $C^g$ is a G-cycle in $D-V(P^g)$. Suppose that for each $i\in[r-1]$, either $u_i$ or $u_iu_{i+1}$ has a partner on $C^g$, and that $u_r$ has a partner on $C^g$, 
		then $D$ has a G-cycle with the vertex set $V(P^g)\cup V(C^g)$ and at least $\ell{}(P^g)+\ell{}(C^g)+1$ arcs.
	\end{lemma}
	
	\pf We proceed by induction on $r$. If $r=1$, then there is nothing to prove, hence assume that $r\geq 2$. Let $xy$ be a partner of $u_1$ or $u_1u_2$ on $C^g$. Choose $i$  maximum such that $u_i\dom y$. Clearly, $P^g[u_1,u_i]$ can be inserted in $C^g$ to give a new G-cycle $L^g$. Moreover, $\ell(L^g)\geq \ell{}(P^g[u_1,u_i])+\ell{}(C^g)+1$. If $i=r$, then we are done. Otherwise, apply induction to the G-path $P^g[u_{i+1},u_r]$ and the G-cycle $L^g$ and note that we may loose an  arc between $u_i$ and $u_{i+1}$ so we cannot guarantee more than one extra arc in total.\qed
	
	We now mention a number of useful known results.
	
	\begin{theorem}\cite{redeiALS7}
		\label{thm:redei}
		Every semicomplete digraph has a Hamiltonian path.
	\end{theorem}
	
	\begin{theorem}\cite{camionCRASP249}
		\label{thm:camion}
		Every strong semicomplete digraph has a Hamiltonian cycle.
	\end{theorem}

	\begin{theorem}\cite{gutinK1988,gutinSJDM6}\label{thm:hampathsmd}
		A semicomplete multipartite digraph $D$ is has a Hamiltonian path  if and only
		if it contains  a 1-path-cycle factor. One can verify whether $D$
		has a Hamiltonian path  and find such a  path in $D$ (if one exists) in time
		$O(n^{2.5})$.
	\end{theorem}

	\begin{theorem}\cite{gutinSJDM6}\label{thm:longpathsmd}
		Let $D$ be a semicomplete multipartite digraph of order $n$.
		\item[(a)] Let $\cal F$ be a 1-path-cycle subdigraph with maximum
		number of vertices in $D$. Then $D$ contains a path $P$ such that
		$V(P)= V({\cal F})$.
		
		\item[(b)] A longest path in $D$ can be constructed in time
		$O(n^3)$.
	\end{theorem}
	
	\begin{theorem}\cite{gutinDM141}\label{thm:extSDhamchar}
		An extended semicomplete digraph is Hamiltonian if and only if it is strong and has a cycle factor.
	\end{theorem}

	\begin{theorem}\label{thm:esdlongc}\cite{gutinthesis}
		Let $D$ be a strong extended semicomplete  digraph and let $\cal
		F$ be a cycle subdigraph of $D$. Then $D$ has a cycle $C$ which
		contains all vertices of $\cal F$. The cycle $C$ can be found in
		time $O(n^3)$. In particular, if $\cal F$ is maximum, then $V(C)=
		V({\cal F})$, i.e., $C$ is a longest cycle of $D$.
	\end{theorem}

	\section{Basic observations about G-cycles and G-paths}\label{Gpathcycle}

\begin{figure}[H]
\begin{center}
\tikzstyle{vertexL}=[circle,draw, top color=gray!30, bottom color=gray!5, minimum size=12pt, scale=1.2, inner sep=0.9pt]
\begin{tikzpicture}[scale=0.3]
  \node (a) at (1,6) [vertexL]{$a$};
  \node (b) at (11,5.3) [vertexL]{$b$};
  \node (c) at (21,6) [vertexL]{$c$};
  \node (x) at (11,1) [vertexL]{$x$};
  \node (y) at (11,10) [vertexL]{$y$};

  \draw[->, line width=0.03cm] (a) to (b);
  \draw[->, line width=0.03cm] (b) to (c);
  \draw[->, line width=0.03cm] (a) to [out=10, in=170]  (c);

  \draw[->, line width=0.03cm] (x) to (a);
  \draw[->, line width=0.03cm] (b) to (x);
  \draw[->, line width=0.03cm] (x) to (c);

  \draw[->, line width=0.03cm] (a) to (y);
  \draw[->, line width=0.03cm] (y) to (b);
  \draw[->, line width=0.03cm] (c) to (y);
\end{tikzpicture} 
\end{center}
\caption{A strong semicomplete multipartite digraph of order $5$ with no $4$-cycle containing vertices from all partite sets $\{a\}, \{b\}, \{c\}, \{x,y\}$.}\label{Ex1}
\end{figure}
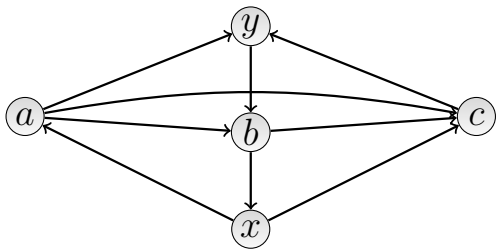

	Recall that the length $\ell{}(C^g)$ of a G-cycle $C^g$ is the number of arcs on $C^g$ and that a longest G-cycle is one with the maximum length.
Consider the \smd{} $D$ with vertex set  $V(D)=\{a,b,c,x,y\}$, arc set $A(D)=\{ab,bc,ac,xa,xc,bx,ay,cy,yb\}$ and partite sets $\{x,y\}$, $\{a\}$, $\{b\}$, $\{c\}$
(see Figure~\ref{Ex1}).
          $D$ has a hamiltonian cycle  $cybxac$ but it has no cycle of length 4 which visits each partite set (every  cycle containing $a$ and $c$ must contain the arcs $cy$ and $xa$ as  $d^+(c)=1$ and  $d^-(a)=1$, implying that the cycle must contain both  $x$ and  $y$). Thus a strong \smd{}
          $D=(V_1,\ldots{},V_c,A)$ may not contain a cycle of length $c$ which visits each partite set.  The following result shows that if we replace 'cycle' by 'G-cycle', then  one can in fact always find a G-cycle of length equal to the number of partite sets and which visits each of these.\\

		\begin{proposition}
			\label{shortGcycle}
			Every strong \smd{} $D=(V_1,\ldots{},V_c,A)$ contains a good G-cycle of length $c$
		\end{proposition}
		\pf
		Consider the digraph $D^c$ that we obtain by contracting each $V_i$, $i\in [c]$ to one vertex and deleting parallel copies of arcs.
		Then $D^c$ is a strong semicomplete digraph on $c$ vertices and hence has a Hamiltonian cycle $C$ by Theorem \ref{thm:camion}.
		Back in $D$ the arcs of $C$ form either a cycle of $D$ or a collection of $r\geq 2$  vertex-disjoint paths $P_{(u_1,v_1)},P_{(u_2,v_2)},\ldots{},P_{(u_r,v_r)}$ which form a good G-cycle of $D$ containing exactly $c$ arcs of $D$. \qed
		
		\begin{proposition}
			\label{shortGpath}
			Every  \smd{} $D=(V_1,\ldots{},V_c,A)$ contains a good G-path of length $c-1$
		\end{proposition}
		\pf Consider again the digraph $D^c$ above. By Theorem \ref{thm:redei} $D^c$ has a Hamiltonian path $P$. Back in $D$ the path $P$ corresponds to a collection of one or more vertex-disjoint paths which form a good G-path containing exactly $c-1$ arcs of $D$.\qed\\

	As we shall emphasize later, finding a longest G-cycle in a \smd{} is not simple. In contrast to this we have the following observation, where a G-cycle is {\bf maximum} if is a longest G-cycle and it contains the maximum number of vertices among all longest G-cycles of $D$.
	
	\begin{lemma}\label{longestG-cycle}
		Let $D$ be a strong \smd{}. Then every maximum G-cycle is a spanning G-cycle.
	\end{lemma}
	\pf 
Let $C^g$ be a maximum G-cycle.  Suppose $C^g$ is not spanning and that $V$ is a partite set such that the  vertex $v\in V$ is not on $C^g$. 

Assume that $v$ has an in-neighbor and an out-neighbor on $C^g$. If $C^g$ contains a vertex $u$ such that $uv,vu\in A(D)$ and all vertices in $V(C^g)-\{u\}$ belong to $V$, then $C^g[u,u^-]vu$ is a G-cycle of length $\ell{}(C^g)$ with more vertices than $C^g$, a contradiction.  So we may assume that there is a path $x_1 v x_2$ in $D$ with $x_1,x_2 \in V(C^g)$. Choose the path $x_1 v x_2$ such that $C^g[x_1,x_2]$ has minimum
order. If $|V(C^g[x_1,x_2])|=2$ then $C^g[x_2,x_1]vx_2$ is a longer G-cycle, a contradiction again.
So, $|V(C^g[x_1,x_2])|>2$ and $x_2^- \in V$. Now $C^g[x_2,x_2^-]vx_2$ is a G-cycle of length at least $\ell{}(C^g)$ with more vertices than $C^g$, a contradiction. 

Therefore, without loss of generality we may assume that $v$ has no out-neighbor on $C^g$. Let $P$ be a $(v,C^g)$-path in $D$, which exists as $D$ is strong, and let $u$ be the last vertex of $P$.
No matter whether $u^- \in V$ or not we observe that
$C^g[u,u^-]P$ is a longer G-cycle
, a contradiction.
So $C^g$  must be spanning, which completes the proof.
\qed

\section{Spanning  G-paths with prescribed end-vertices}

A digraph $D=(V,A)$ is {\bf Hamiltonian-connected} if it contains a spanning $(x,y)$-path for all pairs of distinct vertices $x,y\in V$.

\begin{theorem}\cite{thomassenJCT28}\label{thm:hamconSD}
	Every 4-strong semicomplete digraph is Hamiltonian-connected.
\end{theorem}

 As every strongly connected Hamiltonian-connected digraph has a Hamiltonian cycle and there exists a $k$-strong semicomplete multipartite digraph with no Hamiltonian cycle for every $k\geq1$, see \cite[page 247]{bang2009}, and Figure \ref{noclose} it is easy to see that there can be no $k$ such that every $k$-strong \smd{} is Hamiltonian-connected. When we relax to G-paths we can obtain the following generalization of Theorem \ref{thm:hamconSD}. 

\begin{corollary}
	\label{cor:4stronghamGpath}
	Every 4-strong \smd{} has a spanning $(x,y)$-G-path for every choice of distinct vertices $x,y$.
\end{corollary}
\pf
Let $D=(V_1,\ldots{},V_c,A)$ be 4-strong and let $x,y$ be distinct vertices of $D$. Let $D^*$ be obtained from $D$ by adding all possible arcs inside each $V_i$, that is, $D^*[V_i]$ is a complete digraph for each $i\in [c]$. Then $D^*$ is semicomplete  and  is $4$-strong since already $D$ is 4-strong. Hence it follows from  Theorem \ref{thm:hamconSD} that  $D^*$ has an $(x,y)$-Hamiltonian path $P$. By deleting all arcs of $P$ with both ends in the same partite sets we obtain an $(x,y)$-G-path covering all vertices of $D$. \qed

\begin{theorem}\cite{bangJA13}\label{thm:hamconalg}
	There exists a polynomial algorithm for deciding whether a given semicomplete digraph with distinct vertices $x,y$ has a Hamiltonian $(x,y)$-path.
\end{theorem}

The following corollary of Theorem \ref{thm:hamconalg} is proved in the same way as we proved Corollary \ref{cor:4stronghamGpath}.
\begin{corollary}\label{cor:alghamGpath}
	There exists a polynomial algorithm for deciding whether a given semicomplete multipartite digraph with distinct vertices $x,y$ has a spanning  $(x,y)$-G-path.
\end{corollary}

Note that the construction we used to prove Corollaries \ref{cor:4stronghamGpath} and \ref{cor:alghamGpath} does not always allow us to find a spanning  $(x,y)$-G-path containing the maximum number of arcs.
In fact, the complexity of finding a spanning  
$(x,y)$-G-path with the maximum number of arcs is open, already for semicomplete digraphs (where the path will be a real path).

\begin{conjecture}\cite[Conjecture 7.3.7]{bang2009}
	There exists a polynomial algorithm for finding a longest $(x,y)$-path when the input is a semicomplete digraph $D$ and $x,y$ are distinct vertices of $D$.
\end{conjecture}

\section{Longest G-paths}\label{sec:longestGP}

While the Hamiltonian cycle problem is polynomial
for \smd{}s  \cite{bangJGT29}  (see Theorem \ref{thm:smdhcalg}) there is no easy characterization of \smd{}s with a Hamiltonian cycle. This is sharp contract to the following beautiful result by Gutin.

\begin{theorem}\cite{gutinK1988,gutinSJDM6}\label{thm:longpathSMD}
	The number of vertices of a longest path in a \smd{} $D$ equals the maximum number of vertices that can be covered by a 1-path-cycle subdigraph of $D$. 
\end{theorem}

In fact Gutin showed that the stronger statement that given any 1-path-cycle subdigraph  $P,C_1,\ldots{},C_k$ in a \smd{} one can produce a Hamiltonian path of $D[V(P)\cup{}V(C_1)\cup\cdots\cup{}V(C_k)]$ in polynomial time.\\

The following result generalizes Theorem \ref{thm:longpathsmd} and implies that we can find a longest G-path in any \smd{} in polynomial time.

\begin{theorem}\label{thm:longG-path}
	Let $D=(V_1,\ldots{},V_c,A)$ be a \smd{}. The maximum number of arcs in a G-path of $D$ is equal to the maximum number of arcs contained in a 1-G-path-G-cycle subdigraph. Furthermore, we can find a longest G-path in polynomial time.
\end{theorem}

This follows from the following lemma by induction over the number of G-cycles in the 1-G-path-G-cycle subdigraph.

\begin{lemma}\label{pathcycle-cycle}
	Let $D$ be a \smd{}. Suppose that $P^g=x_1\cdots{}x_s$ is a G-path and $C^g=y_1\cdots{}y_ty_1$ is a G-cycle in $D-V(P^g)$. Then in polynomial time we can produce a  G-path $Q^g$ in $D$ such that $V(Q^g)=V(P^g)\cup V(C^g)$ and $\ell{}(Q^g)\geq \ell{}(P^g)+\ell{}(C^g)$.
\end{lemma}
\pf Suppose to the contrary that there is no such G-path in $D$. First observe that there is no vertex $y_i\in V(C^g)$ which dominates $x_1$, that is, $N^-(x_1)\cap V(C^g)=\emptyset$  otherwise, $C^g[y_{i+1},y_i]P^g$ is the wanted G-path. Analogously, we have $N^+(x_s)\cap V(C^g)=\emptyset$.\\

If there are indices $a\in [2,s]$ and $b\in [t]$ such that $y_b\dom x_a$ and $x_{a-1}\dom y_{b+1}$, then $Q^g=P^g[x_1,x_{a-1}]C^g[y_{b+1},y_b]P^g[x_a,x_s]$ is a G-path with $V(Q^g)=V(P^g)\cup V(C^g)$ and $\ell{}(Q^g)\geq \ell{}(P^g)+\ell{}(C^g)$. If instead we have $y_b\dom x_a$ and $x_{a-1},y_{b+1}$ belong to the same partite set, then the G-path $Q^g$ above still satisfies $\ell{}(Q^g)\geq \ell{}(P^g)+\ell{}(C^g)$ unless $x_{a-1}x_a,y_by_{b+1}\in A(D)$. Hence we may assume that the following holds.

\begin{claim}\label{clm1}
	Suppose that $y_bx_a\in A(D)$.  Then  $a>1$ and either  $y_{b+1}\dom x_{a-1}$ or $y_{b+1}$ and $x_{a-1}$ belong to the same partite set 
	and $x_{a-1}x_a,y_by_{b+1}\in A(D)$.
\end{claim}

Suppose first that there are no indices $i,j$ such that both $y_{j+1}x_i$ and $y_jx_{i+1}$ belong to $A(D)$. If there is no arc between $C^g$ and $x_s$, then all vertices of $C^g$ and $x_s$ belong to the same partite set, in particular $C^g$ has no arcs. Thus  $P^gC^g[y_1,y_t]$ is the wanted G-path. As $N^+(x_s)\cap V(C^g)=\emptyset$, we may assume that $y_r$ dominates $x_s$ for some $y_r\in V(C^g)$. Now it follows from the assumption above and Claim \ref{clm1} that $x_{s-1}$ and $y_{r+1}$ belong to the same partite set $V$ and that both $x_{s-1}x_s$ and $y_ry_{r+1}$ belong to $A(D)$. As $N^+(x_s)\cap V(C^g)=\emptyset$ again, we have $y_{r+1}$ dominates $x_s$. Thus, again the assumption above  and Claim \ref{clm1} implies that we also have $y_{r+2}\in V$. Now $Q^g=P^g[x_1,x_{s-1}]C^g[y_{r+2},y_{r+1}]x_s$ is the wanted G-path.

It remains to consider the case that there exist $i$ and $j$ such that $\{y_{j+1}x_i, y_jx_{i+1}\}\subseteq A(D)$. Choose $i$ as small as possible. Clearly, $i>1$ as $N^-(x_1)\cap V(C^g)=\emptyset$. Then $y_{j+2}x_{i-1}\notin A(D)$ follows by the minimality of $i$ and then  Claim \ref{clm1} implies that  $y_{j+2}$ and $x_{i-1}$ belong to the same partite set $V'$ and that $x_{i-1}x_i,y_{j+1}y_{j+2}\in A(D)$. Then using that $y_{j+2}$ and $x_{i-1}$ belong to the same partite set, we conclude that there is an arc between $x_{i-1}$ and $y_{j+1}$, and an arc between $x_{i}$ and $y_{j+2}$.

Suppose that $y_{j+2}$ dominates $x_{i}$.  Then it follows from the minimality of $i$ and Claim \ref{clm1} that $x_{i-1}$ and $y_{j+3}$ belong to the same partite set and that $y_{j+2}y_{j+3}$ belong to $A(D)$, which is impossible as $x_{i-1},y_{j+2},y_{j+3}$ all belong to the same partite set. Hence we have $x_i\dom y_{j+2}$. We claim that we also have $x_{i-1}\dom y_{j+1}$. 
This is clear if  $i=2$ as $N^-(x_1)\cap V(C^g)=\emptyset$ so suppose that  $i>2$ and $y_{j+1}x_{i-1}\in A(D)$. Then the minimality of $i$ and Claim \ref{clm1} implies that $x_{i-2}$ and $y_{j+2}$ belong to the same partite set  and $x_{i-2}x_{i-1}\in A(D)$, which is impossible as  $x_{i-1},x_{i-2},y_{j+2}$ all belong to the same partite set. Thus we must have $x_{i-1}\dom y_{j+1}$ and  now
$Q^g=P^g[x_1,x_{i-1}]y_{j+1}x_iC^g[y_{j+2},y_j]P^g[x_{i+1},x_s]$ is the wanted G-path (we gain 4 arcs between $P^g$ and $C^g$ and loose at most 4 arcs). \qed
\begin{figure}[H]
	\begin{center}
		\tikzstyle{vertexL}=[circle,draw, minimum size=18pt, scale=0.85, inner sep=0.5pt]
		\tikzstyle{vertexS}=[circle,draw, minimum size=18pt, scale=0.6, inner sep=0.5pt]
		\tikzstyle{vertexSx}=[circle,draw, top color=black!30, bottom color=black!5, minimum size=18pt, scale=0.6, inner sep=0.5pt]
		\begin{tikzpicture}[scale=0.4]
		\node at (0,10) {  };
		\node (x1) at (1,8) [vertexL]{$x_1$};
		\node (x2) at (4,8) [vertexL]{$x_2$};
		\node (xim1) at (8,8) [vertexSx]{$x_{i-1}$};
		\node (xi) at (11,8) [vertexL]{$x_i$};
		\node (xip1) at (14,8) [vertexS]{$x_{i+1}$};
		\node (xsm1) at (18,8) [vertexS]{$x_{s-1}$}; 
		\node (xs) at (21,8) [vertexL]{$x_s$};
		
		\node (y1) at (21,2) [vertexL]{$y_1$};
		\node (y2) at (18,2) [vertexL]{$y_2$};
		\node (yj) at (14,2) [vertexL]{$y_j$};
		\node (yjp1) at (11,2) [vertexS]{$y_{j+1}$};
		\node (yjp2) at (8,2) [vertexSx]{$y_{j+2}$};
		\node (ytm1) at (4,2) [vertexS]{$y_{t-1}$};
		\node (yt) at (1,2) [vertexL]{$y_t$};
		
		\draw[->, line width=0.06cm] (x1) to (x2);
		\draw[line width=0.06cm] (x2) to (5.2,8);   \node at (6,8) {{\small $\cdots$}};   \draw[->, line width=0.06cm] (6.8,8) to (xim1);
		\draw[->, line width=0.01cm] (xim1) to (xi);
		\draw[->, line width=0.01cm] (xi) to (xip1);
		\draw[line width=0.06cm] (xip1) to (15.2,8);   \node at (16,8) {{\small $\cdots$}};   \draw[->, line width=0.06cm] (16.8,8) to (xsm1);
		\draw[->, line width=0.06cm] (xsm1) to (xs);
		
		\draw[->, line width=0.06cm] (y1) to (y2);
		\draw[line width=0.04cm] (y2) to (16.8,2);   \node at (16,2) {{\small $\cdots$}};   \draw[->, line width=0.06cm] (15.2,2) to (yj);
		\draw[->, line width=0.01cm] (yj) to (yjp1);
		\draw[->, line width=0.01cm] (yjp1) to (yjp2);
		\draw[line width=0.06cm] (yjp2) to (6.8,2);   \node at (6,2) {{\small $\cdots$}};   \draw[->, line width=0.06cm] (5.2,2) to (ytm1);
		\draw[->, line width=0.06cm] (ytm1) to (yt);
		
		\draw[->, line width=0.06cm] (yjp1) to (xi);
		\draw[->, line width=0.06cm] (yj) to (xip1);
		\draw[->, line width=0.06cm] (xim1) to (yjp1);
		\draw[->, line width=0.06cm] (xi) to (yjp2);
		\draw[line width=0.06cm] (yt)  to [out=320, in=180]  (4,0.3);
		\draw[line width=0.06cm] (4,0.3) to  (18,0.3);
		\draw[->, line width=0.06cm] (18,0.3)  to [out=0, in=220]  (y1);
		\end{tikzpicture}
	\end{center}
	\caption{An illustration of the final G-cycle in Lemma~\ref{pathcycle-cycle}, where the darker vertices belong to the same partite set, and the thin arcs are not included in the G-cycle.}\label{Thm52}
\end{figure}

\section{Good G-cycles in extended semicomplete digraphs}

In this  section $D=S[I_{n_1},\ldots{},I_{n_s}]$ denotes an extended semicomplete digraph. This means that the $i$th vertex of a semicomplete digraph $S$ was replaced by an independent set $I_{n_i}$ of $n_i\geq 1$ vertices.
We call a pair of distinct vertices from the same partite set of $D$ {\bf similar} (they have exactly the same in- and out-neighbours as $D$ is extended semicomplete).

The following result generalizes Theorem \ref{thm:extSDhamchar} and together with Theorem \ref{thm:findlongestgoodGC} it implies that the longest G-cycle problem can be solved in polynomial time when the input is an extended semicomplete digraph.

\begin{theorem}\label{thm:longGCextSD}
	An extended semicomplete digraph $D$ has a spanning G-cycle if and only if it is strong and has a G-cycle factor. Such a cycle can be produced on polynomial time from a given G-cycle factor. Furthermore, the maximum number of arcs in spanning G-cycle is equal to the maximum number of arcs contained in a $G$-cycle factor of $D$.
\end{theorem}

To prove this result we first establish two  lemmas both  of which generalize results of \cite{gutinDM141} for extended semicomplete digraphs.
The proofs of these lemmas are similar to the original proofs, but we need to handle the  cases when there are  consecutive vertices from the same partite set in one or more of the cycles.

\begin{lemma}\label{lem:similarmerge}
	If $C_1^g$ and $C_2^g$ are G-cycles which both contain a vertex from the same partite set $V$ of an extended semicomplete digraph $D$,
	then $D$ has a G-cycle $C^g$ with $V(C^g)=V(C^g_1)\cup V(C^g_2)$ and  $\ell{}(C^g)=\ell{}(C^g_1)+\ell{}(C^g_2)$.
\end{lemma}

\pf Let $C_1^g=u_1\cdots{}u_pu_1$ and let $C_2^g=v_1\cdots{}v_qv_1$. If $\ell{}(C^g_i)=0$ (it has no arcs), then all vertices of $C^g_i$ belong to  $V$ and  the claim is trivial as we can replace any vertex $v\in V$ on $C^g_{3-i}$ by the vertices $V(C^g_i)\cup\{v\}$. So we can assume that $\ell{}(C^g_i)\geq 1$ for $i=1,2$. 
Let $r\in [p]$ and $s\in [q]$ be chosen such that $u_{r}$ and $v_{s}$ belong to the same partite set $V$ and $u_{r-1}\dom u_{r}$, $v_{s-1}\dom v_{s}$ both hold (this is possible as not all vertices of $C^g_i$ belongs to $V$).
Then $v_{s-1}\dom u_{r}$ and $u_{r-1}\dom v_{s}$ as $u_{r}, v_{s}$ are similar and $D$ is extended semicomplete. Now it follows that $C_1^g[u_r,u_{r-1}]C_2^g[v_s,v_{s-1}]u_r$ is G-cycle $C^g$ with $V(C^g)=V(C^g_1)\cup V(C^g_2)$ and 
$\ell{}(C^g)=\ell{}(C^g_1)+\ell{}(C^g_2)$. \qed

\begin{lemma}\label{lem:mergetwoGC}
	If $C_1^g$ and $C_2^g$ are G-cycles such that the extended semicomplete digraph $D$ contains some arc from $C^g_i$ to $C^g_{3-i}$ for $i=1,2$, then
	$D$ has a G-cycle $C^g$ with $V(C^g)=V(C^g_1)\cup V(C^g_2)$ and such that $\ell{}(C^g)\geq \ell{}(C^g_1)+\ell{}(C^g_2)$.
\end{lemma}

\pf Let $C_1^g=u_1\cdots{}u_pu_1$ and let $C_2^g=v_1\cdots{}v_qv_1$. By Lemma \ref{lem:similarmerge} we may assume that there is no pair of vertices $u_{i},v_{j}$ which are similar. This implies that every vertex on $C^g_1$ is adjacent to every vertex of $C^g_2$.  Suppose first that there is an $i\in [2]$ such that $C^g_i$ has a vertex which has an arc to every vertex of $C^g_{3-i}$. Without loss of generality $i=1$ and let $u_a$ be such a vertex. By the assumption there is also a vertex of $C^g_2$ which sends an arc to $V(C^g_1)$. Hence we may assume that $u_{a}$ is chosen such that there is an arc $v_{b}u_{a+1}$ in $D$. As $u_{a}v_{b+1}$ is an arc (by the choice of $u_{a}$) we get the G-cycle $C^g=u_{a}C^g_2[v_{b+1},v_{b}]C^g_1[u_{a+1},u_{a}]$. It is easy to check that $\ell{}(C^g)\geq \ell{}(C^g_1)+\ell{}(C^g_2)$.

Thus we may assume that every vertex of $C^g_i$ has arcs in  both directions to $V(C^g_{3-i})$ for all $i\in [2]$. 
Now we claim that for every vertex $v_{r}$ of $C^g_2$ there exist consecutive vertices
$u_{s},u_{s+1}$ such that $u_{s}v_{r}$ and $v_{r}u_{s+1}$ are both arcs of $A(D)$.
Suppose not and let $u_{i}$, $u_{i+d}$ be chosen such that $u_{i}v_{r}$ and
$v_{r}u_{i+d}$ are both arcs of $A(D)$ and $d\geq 1$ is as small as possible. By the assumption, $d>1$, but then $u_{i+1}$ and $v_{r}$ must be non-adjacent and hence similar, contradicting the assumption above. Thus we have shown that every vertex of $C^g_2$ has a partner on $C^g_1$ and now it follows by Lemma \ref{partner-spanning} that $D$ has a G-cycle 
$C^g$ with $V(C^g)=V(C^g_1)\cup V(C^g_2)$ and such that $\ell{}(C^g)\geq \ell{}(C^g_1)+\ell{}(C^g_2)$.\qed\\

\noindent{}{\bf Proof of Theorem \ref{thm:longGCextSD}:}
Let us first show that if an extended semicomplete digraph $D$ with more than one partite set has a spanning G-cycle, then it must be strong.
Let $C^g$ be a spanning G-cycle of $D$ and 
let $u$ and $v$ be arbitrary vertices in $D$. If $C^g[u,v]$ only contain vertices from one partite set in $D$, then swap $u$ and $v$ on $C^g$ (which is
possible as $u$ and $v$ are similar).
Now let $P=C^g[u,v]$, and note that $P$ contains vertices from at least two partite sets.
 Everytime there are two consecutive vertices $z$ and $z^+$ 
on $P$ from the same partite set we consider the path $P[u,z] P[z^{++},v]$ (if $z^+ \not= v$) or $P[u,z^-] P[z^{+},v]$ (if $z \not=u$) instead of $P$.
Continuing this process until there are no more consecutive vertices on $P$ from the same partite set, gives us a $(u,v)$-path in $D$.
As $u$ and $v$ were arbitrary, $D$ is therefore strong.

%
%

Let $F=C^g_1\cup\cdots{}\cup C^g_p$, $p\geq 1$ be a G-cycle factor of $D$. By Lemmas \ref{lem:similarmerge} and \ref{lem:mergetwoGC} we can assume that every two vertices on distinct G-cycles $C^g_i,C^g_j$  of $F$ are adjacent and either 
 $C^g_i\compdom C^g_j$ or $C^g_j\compdom C^g_i$. 
This means that if we contract each G-cycle $C^g_i$ to a vertex $x_i$ for $i\in [p]$ and delete loops and parallel arcs resulting from this operation we are left with a tournament $T$ on $p$ vertices. Further, as $D$ is strong so is $T$. By Theorem \ref{thm:camion}, $T$ has a Hamiltonian cycle and by changing the labels of the G-cycles in $F$ we can assume that $x_1x_2\cdots{}x_px_1$ is a Hamiltonian cycle of $T$. This means that back in $D$ we have $C^g_1\compdom C^g_2\compdom\cdots\compdom{}C^g_p\compdom C^g_1$. Let $C_i^g=v_{i,1}v_{1,2}\cdots{}v_{i,q_i}v_{i,1}$ for all $i\in[p]$. Then $C^g=C_1^g[v_{1,1},v_{1,q_1}]C_2^g[v_{2,1},v_{2,q_2}]\cdots{}C_p^g[v_{p,1},v_{p,q_p}]v_{1,1}$ is a spanning G-cycle of $D$ and $\ell{}(C^g)\geq \ell{}(C^g_1)+\cdots{}+\ell{}(C^g_p)$ as we add $p$ new arcs and loose at most one arc per cycle. \qed

\begin{theorem}\label{thm:findlongestgoodGC}
	There exists a polynomial algorithm for checking whether a given strongly connected extended semicomplete digraph has a spanning G-cycle and if so, then find such a cycle containing the maximum number of arcs.
\end{theorem}

\pf Since the proof above provides an algorithm for producing a spanning cycle with at least $r$ arcs
from a given G-cycle factor with $r$ arcs, it suffices to recall that by Lemma \ref{lem:maxnoarcsCF} we can find a G-cycle factor containing the maximum number of arcs of $D$. \qed

\section{G-cycle factors with minimum number of G-cycles}

In this section we study the structure of G-cycle factors in general \smd{}s and obtain several results which generalize previous results from \cite{bangDM161,yeoJGT24} on cycle factors in \smd{}s. We point our here that all proofs in this section are constructive and  can be converted to polynomial algorithms for obtaining the desired objects.

Recall that a  \textbf{G-cycle factor} of a \smd{} $D$ is a collection of vertex-disjoint G-cycles covering all vertices of $D$.
A G-cycle factor $F=C^g_1\cup\cdots{}\cup C^g_t$ of a \smd{} $D$ is {\bf irreducible} if $D$ has no G-cycle factor
$F'=\tilde{C}^g_1\cup\cdots{}\cup \tilde{C}^g_r$ such that $|A(F')|\geq |A(F)|$ and $r<t$ and for every $j\in [t]$ we have $V(C^g_j)\subseteq V(\tilde{C}^g_i)$ for some $i\in [r]$.
For a G-cycle $C^g$ and a vertex $v\notin V(C^g)$, we say that $v$ is \textbf{out-singular} (resp., \textbf{in-singular}) with respect to $C^g$ if $v\compdom C^g$ (resp., $C^g\compdom v$), $v$ is \textbf{singular} with respect to $C^g$ if it is either out-singular or in-singular with respect to $C^g$.

Let $D$ be a digraph with two vertex-disjoint G-cycles $C^g_1$ and $C^g_2$. We use the notation $C^g_1\simeq> C^g_2$ to denote that $C^g_1$ contains singular vertices with respect to $C^g_2$ and they all are out-singular, and $C^g_2$ contains singular vertices with respect to $C^g_1$ and they all are in-singular.

A G-cycle factor $F=C^g_1\cup\cdots{}\cup C^g_t$ is called \textbf{feasible} if it has no pair of G-cycles $C^g_i,C^g_j (i\neq j)$ such that
$C^g_i\simeq> C^g_j$ or $C^g_j\simeq> C^g_i$. 

\bigskip


\begin{lemma}\label{nonsingular}
	Let $C^g\cup Q^g$ be a G-cycle factor in a \smd{} $D$. Let $C^g=u_1u_2\cdots{}u_pu_1, p\geq 2$ and suppose that $D$ has no spanning G-cycle with at least $\ell{}(Q^g)+\ell{}(C^g)$ arcs. Then for every pair  $v,v^+$ of consecutive vertices on 
	$Q^g$, if $v$ is not singular with respect to $C^g$, then either $v$ has a partner on $C^g$ or, $vv^+$ is an arc and it has a partner on $C^g$.
\end{lemma}
\pf Suppose that $v$ has no partner on $C^g$. Since $v$ is non-singular, there exist $u_s,u_t\in V(C^g)$ such that $u_s\dom v\dom u_{t}$. Choose a pair of vertices $u_i,u_{i+d}$ with $u_{i}\dom v\dom u_{i+d}$ such that $d\geq1$ is minimum among all such pairs. Note that $2\leq d\leq p$ as $v$ has no partner on $C^g$ and if $d=p$ then $u_i=u_{i+d}$. Moreover, by the minimality of $d$, all vertices in $C^g[u_{i+1},u_{i+d-1}]$ and $v$ belong to the same partite set $V$.


If  $v^+$ also belongs to $V$, then $Q^g[v^+,v]C^g[u_{i+d},u_{i+d-1}]v^+$ is a spanning G-cycle of $D$ with $\ell{}(Q^g)+\ell{}(C^g)$ arcs, contradicting the assumption.  So $v^+\notin V$ and then $vv^+$ is an arc of $Q^g$. Moreover, there is an arc between $u_{i+1}$ and $v^+$ as $u_{i+1}\in V$. If $u_{i+1}\dom v^+$, then $Q^g[v^+,v]C^g[u_{i+2},u_{i+1}]v^+$ is a spanning G-cycle containing at least $\ell{}(Q^g)+\ell{}(C^g)$ arcs, a contradiction.  Hence we must have $v^+\dom u_{i+1}$, showing that $u_iu_{i+1}$ is a partner of the arc $vv^+$.\qed

\begin{lemma}\cite{bangDM161}\label{jbj1996}
	Let $D$ be a \smd{} containing a (real) cycle factor $C_1\cup C_2$ such that $C_i$ has no singular vertices with respect to $C_{3-i}$ for each $i\in[2]$, then $D$ is Hamiltonian. 
\end{lemma}
\begin{lemma}
	Let $D$ be a \smd{} containing a G-cycle factor $C^g_1\cup C^g_2$ such that $C^g_i$ has no singular vertices with respect to $C^g_{3-i}$ for each $i\in[2]$, then $D$ has a spanning G-cycle  with at least $\ell{}(C^g_1)+\ell{}(C^g_2)$ arcs.
\end{lemma}
\pf Suppose to the contrary that $D$ has no such a cycle. By  Lemma \ref{nonsingular}, 
for each $i\in [2]$ and for each pair of consecutive vertices $w,w^+$ of  $C^g_i$, either $w$ has a partner on $C^g_{3-i}$ or $ww^+$ is an arc. 

If there exists a vertex $w$ on $C^g_j$ with some $j\in[2]$ such that $w$ has a partner on $C^g_{3-j}$, then by Lemma  \ref{partner-spanning} we can merge $C^g_j[w^+,w]$ and $C^g_{3-j}$ into a spanning G-cycle with length at least $\ell{}(C^g_j[w^+,w])+\ell{}(C^g_{3-j})+1\geq \ell{}(C^g_1)+\ell{}(C^g_2)$, which contradicts our assumption.

Hence  we may assume that $C_1^g$ and $C_2^g$ are real cycles by Lemma \ref{nonsingular} and Lemma \ref{jbj1996} shows that $D$ is Hamiltonian, which contradicts our assumption again and this completes the proof.\qed



\begin{lemma}
	Let $D$ be a \smd{} containing a G-cycle factor $C^g_1\cup C^g_2$ such that only $C^g_1$ has singular vertices with respect to $C^g_2$, then $D$ has a spanning G-cycle with at least $\ell{}(C^g_1)+\ell{}(C^g_2)$ arcs.
\end{lemma}
\pf By reversing all arcs if necessary, we may assume that $C^g_1$ has an out-singular vertex $x$ and $C^g_2$ has no singular vertices. Since $C^g_2$ has no singular vertices, $C^g_1$ has a  vertex which is not out-singular. Choose $x$ such that $x$ is out-singular but its successor $x^+$ is not an out-singular vertex and assume that $y\in C^g_2$ dominates $x^+$.

As $C^g_2$ has no singular vertices, by  Lemma \ref{nonsingular}, for any two consecutive vertices $u_i,u_{i+1}$ on $C^g_2$, either $u_i$ or $u_iu_{i+1}$ has a partner on $C^g_1$. If $y$ has a partner on $C^g_1$, then the lemma follows by applying Lemma \ref{partner-spanning} with $C^g_2[y^+,y]$ and $C^g_1$. So we may assume that $x$ does not dominate $y$, which implies that $x$ and $y$ belong to the same partite set. Therefore,  $xC^g_2[y^+,y]C^g_1[x^+,x]$ is a spanning G-cycle with at least $\ell{}(C^g_1)+\ell{}(C^g_2)$ arcs. Note that this holds whether or not $x$ and $y^+$ belong to the same partite set or $x$ dominates $y^+$.\qed

\begin{lemma}\label{goodsingular}
	Let $D$ be a \smd{} containing a G-cycle factor $C^g_1\cup C^g_2$. Suppose that both $C^g_1$ and $C^g_2$ have singular vertices, but $(C^g_1,C^g_2)$ is feasible, then $D$ has a spanning G-cycle with at least
	$\ell{}(C^g_1)+\ell{}(C^g_2)$ arcs.
\end{lemma}

\pf Assume w.l.o.g that $C^g_1$ has an out-singular vertex $u$ with respect to $C^g_{2}$. If $C^g_2$ also has an out-singular vertex $v$ with respect to $C^g_{1}$, then $u$ and $v$ are non-adjacent, that is, they belong to the same partite set $V$. Furthermore,  we have $u\dom v^+$ or $v^+\in V$ and, $v\dom u^+$ or $u^+\in V$. Now it follows that the spanning G-cycle $C^g_1[u^+,u]C^g_2[v^+,v]u^+$
has at least $\ell{}(C^g_1)+\ell{}(C^g_2)$ arcs.

Thus we may assume that $C^g_2$ has no out-singular vertex and  then it has an in-singular vertex $v$.
Since $C^g_1\cup C^g_2$ is a feasible G-cycle-factor, $C^g_1$ has both in- and out-singular vertices. Let $w$ be an in-singular vertex on $C^g_1$. By the same argument as above, we have that $w$ and $v$ belong to the same partite set, moreover, either $w^-\dom v$ or $w^-$ and $v$ belong to the same partite set and either $v^-\dom w$ or $v^-$ and $w$ belong to the same partite set and then $C_1^g[w,w^-]C_2^g[v,v^-]w$ is the wanted spanning G-cycle of $D$.\qed

\begin{corollary}\label{2cycleorder}
	Let $D$ be a \smd{} and let $F=C^g_1\cup\cdots{}\cup C^g_k$ be an irreducible G-cycle factor. Then for any two distinct G-cycles $C^g_i,C^g_j$, we have either $C^g_i\simeq> C^g_j$ or $C^g_j\simeq> C^g_i$ (but not both).
\end{corollary}

\begin{lemma}\label{feasible}
	Let $D$ be a \smd{} and let $C^g_1$ and $C^g_2$ be a pair of vertex-disjoint G-cycles in $D$ such that $C^g_1	\simeq> C^g_2$ and $C^g_1\ncompdom C^g_2$. 
Assume that there is no G-cycle in $D$ with vertex set $V(C^g_1)\cup V(C^g_2)$ and  at least $\ell{}(C^g_1)+\ell{}(C^g_2)$ arcs. 
Then there exists a unique partite set $V$ of $D$ such that for every $(V(C^g_2),V(C^g_1))$-G-path $P^g$ starting at a vertex $u$ and 
terminating at vertex $v$, either $\{u_{C^g_2}^+,v_{C^g_1}^-\}\subseteq V$ or there exists a G-cycle $\tilde{C^g}$ in $D$, with $V(\tilde{C^g})=V(C^g_1)\cup V(C^g_2)\cup V(P^g)$ and $\ell{}(\tilde{C^g})\geq \ell{}(C^g_1)+\ell{}(C^g_2)+\ell{}(P^g)$-1. 
\end{lemma}

\pf  Since $C^g_1	\simeq> C^g_2$ and $C^g_1\ncompdom C^g_2$, there is a vertex $x\in V(C^g_1)$ with $x\compdom C^g_2$ and $x^+\ncompdom C^g_2$. Let $V$ be the partite set containing the vertex $x$ and let $y$ be a vertex such that $y^-\dom x^+$. Then we must have $y\in V$ and  $xx^+,y^-y\in A(D)$ as otherwise, $C^g_2[y,y^-]C^g_1[x^+,x]y$ is a G-cycle in $D$ with vertex set $V(C^g_1)\cup V(C^g_2)$ with at least $\ell{}(C^g_1)+\ell{}(C^g_2)$ arcs, which contradicts our assumption. Now we conclude that $y^-\notin V$ and $x\dom y^-$ and then $xx^+$ is a partner of $y^-$.

\textbf{Case 1. $C^g_1\compdom y$.}

Recall that the G-path $P^g$ starts in $u\in V(C^g_2)$ and ends in $v\in V(C^g_1)$. For simplicity, we write $u^+$ and $v^-$ rather than $u_{C^g_2}^+$ and $v_{C^g_1}^-$. If $\{u^+,v^-\}\cap V=\emptyset$, then $C_1^g[x^+,v^-]C_2^g[y,u]P^gC_1^g[v,x]C_2^g[u^+,y^-]x^+$ is the wanted G-cycle (we loose at most 4 arcs of the cycles and gain 3 arcs between the cycles plus all arcs of $P^g$). So at least one of $u^+,v^-$ must belong to $V$. If  both belong to $V$, we are done, so we may assume that exactly one of $u^+$ and $v^-$ belongs to $V$ and then we can assume that   $u^+\dom v^-$ as otherwise we can take $\tilde{C}^g$ to be the  G-cycle $C^g_2[u^+,u]P^gC^g_1[v,v^-]u^+$. 

Suppose that $v^{--}$ and $u^{++}$ belong to the same partite set $V'$. If $V'=V$, then the G-cycle $v^{--}C^g_2[u^{++},u^+]C^g_1[v^-,v^{--}]$ has  $\ell{}(C^g_1)+\ell{}(C^g_2)$ arcs, contradicting the assumption of the lemma (here we used that either $u^+\in V$ or $v^-\in V$ so we loose only one arc of the cycles and gain the arc $u^{+}v^{-}$.) So we can assume that $V'\neq V$ and 
then  $C^g_1[x^+,v^{--}]C^g_2[y,u^+]C^g_1[v^-,x]C^g_2[u^{++},y^-]x^+$  is a G-cycle in $D$ with vertex set $V(C^g_1)\cup V(C^g_2)$ with at least $\ell{}(C^g_1)+\ell{}(C^g_2)$ arcs, which contradicts our assumption.  Therefore, we may assume that $v^{--}$ and $u^{++}$ are adjacent and then the assumption of the lemma implies that we must have $u^{++}\dom v^{--}$ and one of $u^{++}$ and $v^{--}$ does not belong to $V$. Continuing this process we obtain that $u^+\dom v^-$, $u^{++}\dom v^{--}, \ldots,$ which contradicts the fact that $C^g_1$ has the  vertex $x$ which dominates $C^g_2$.

\textbf{Case 2. $C_1^g\ncompdom y$.} 

Label the vertices in $C^g_2$ such that $C^g_2=y_1y_2\cdots y_my_1$, where $y_1=y$. Recall that $x$ and $y$ both belong to the partite set $V$. Now define the statements $\alpha_K$ and $\beta_K$ as follows:

$\alpha_K:$ The vertex $y_k\in V$ and $C^g_1\ncompdom y_k$ for every $k=\{1,3,5,\ldots,K\}$.

$\beta_K:$  $y_ky_{k+1}$ has a partner on $C^g_1[x^+,x]$ for every $k=\{1,3,5,\ldots,K\}$.
\medskip

Note that it suffices to show that $\alpha_K$ and $\beta_K$ are true for every odd $K$ with $1\leq K<m$. In fact, if this statement holds, by the fact that $xx^+$ is a partner of $y_m=y^-$ (as we established in the beginning of the proof) and by 
	applying Lemma  \ref{partner-spanning} to the G-path $C^g_2[y_1,y_m]$ and G-cycle $C^g_1$, there is a G-cycle in $D$ with vertex set $V(C^g_1)\cup V(C^g_2)$ with at least $\ell{}(C^g_1)+\ell{}(C^g_2)$ arcs, which contradicts the  assumption in the lemma.

Clearly $\alpha_1$ holds as we are in Case 2 and we saw above that $y\in V$. Let us first show that $\alpha_1$ implies $\beta_1$. Applying Lemma \ref{nonsingular} with $v=y_1$, either $y_1y_2$ has a partner on $C^g_1$, that is, $\beta_1$ holds, or $y_1$ has a partner $ww^+$ on $C^g_1$. In the later case the G-cycle $C^g_1[w^+,x]C^g_2[y_2,y_m]C^g_1[x^+,w]y_1w^+$ has at least $\ell{}(C^g_1)+\ell{}(C^g_2)$ arcs, contradicting the assumption of the lemma (if $y_2\in V$, then we gain 3 arcs and loose 3). Thus $\beta_1$ holds.

First observe that  when $\beta_{K-2}$ holds, then $y_K\in V$ must also hold as otherwise $x\dom y_K$ and then $C^g_1[x^+,x]C^g_2[y_{K},y_m]x^+$ is a G-cycle and by $\beta_{K-2}$ and Lemma \ref{partner-spanning} we can insert the subpath $C^g_2[y_1,y_{K-1}]$ into this G-cycle  and obtain a G-cycle with vertex set $V(C^g_1)\cup V(C^g_2)$ with at least $\ell{}(C^g_1)+\ell{}(C^g_2)$ arcs, which contradicts our assumption. Hence $\beta_{K-2}$ implies that $y_K\in V$ holds.

Next we show that $\alpha_K$ and $\beta_{K-2}$ imply $\beta_K$. Suppose first  that $y_K$ has no partner. Then by $\alpha_K$, $C^g_1\ncompdom y_K$ and the fact that $C^g_1\simeq> C^g_2$, the vertex $y_K$ is not singular with respect to $C^g_1$. Now  Lemma \ref{nonsingular} implies that $y_Ky_{K+1}$ has a partner, showing that $\beta_K$ holds.\\

So $y_K$ has a partner on $C^g_1$. 
Now by  applying Lemma \ref{partner-spanning} to the  G-path $C^g_2[y_1,y_K]$ and the G-cycle $C^g_1[x^+,x]C^g_2[y_{K+1},y_m]x^+$,  there is a G-cycle in $D$ with vertex set $V(C^g_1)\cup V(C^g_2)$ with at least
$\ell{}(C^g_1)+\ell{}(C^g_2)$ arcs, which contradicts our assumption. Note that  if $y_{K+1}\in V$, then because we also have $y_K\in V$, we do not loose an arc between $y_K$ and $y_{K+1}$ when we make the new G-cycle above.

The final step is to show that $\alpha_{K-2}$ and $\beta_{K-2}$ imply $\alpha_K$. We have already seen that $y_K\in V$ must hold by $\beta_{K-2}$. It follows from  $\alpha_{K-2}$ and $C^g_1	\simeq> C^g_2$ that the vertex $y_{K-2}$ is not singular. As $\beta_{K-2}$ holds, there exist $u_i$ and $u_{i+1}$ on $C_1^g[x^+,x]$ such that $u_i\dom y_{K-2}$ and $y_{K-1}\dom u_{i+1}$. In particular $u_i\not\in V$ as $y_{K-2}\in V$. If $C^g_1\compdom y_K$, then $u_i$ dominates $y_K$  and now $C^g_1[u_{i+1},u_i]C^g_2[y_K,y_{K-1}]u_{i+1}$ is a G-cycle in $D$ with vertex set $V(C^g_1)\cup V(C^g_2)$ with at least $\ell{}(C^g_1)+\ell{}(C^g_2)$ arcs, which contradicts our assumption. Thus $\alpha_K$ holds and this completes the proof.
\qed\\

In the next three results  $D$ is a \smd{}  with a pair of vertex-disjoint G-cycles $C^g_1,C^g_2$  such that
	\begin{equation}\label{C*}
	\begin{aligned}
		C^g_1\simeq> C^g_2 &\mbox{ and  there is no G-cycle in } D \mbox{ with vertex set }\\V(C^g_1)\cup &V(C^g_2)\mbox{ and at least } \ell{}(C^g_1)+\ell{}(C^g_2)\mbox{ arcs}. 
	\end{aligned}\tag{$C^{\ast}$}
\end{equation}

\begin{lemma}\label{relationship}
	Suppose that (\ref{C*}) holds. For every arc $uv$ from $C^g_2$ to $C^g_1$, the following statements hold.
	
	(i) $u^+,v^-$ all belong to some partite set $V$. 
	
	(ii) $v^-v,uu^+\in A(D)$, which implies that $u,v\notin V$.
	
	(iii) $v^-\dom u$ and $v\dom u^+$.
\end{lemma}
\pf 
(i) This follows from Lemma \ref{feasible} by taking the G-path $P^g$ to be the arc $uv$.

(ii) If one of $v^-v$ and $uu^+$ does not belong to $A(D)$, then $v^-C^g_2[u^+,u]C^g_1[v,v^-]$ is a G-cycle in $D$ with vertex set $V(C^g_1)\cup V(C^g_2)$ with at least $\ell{}(C^g_1)+\ell{}(C^g_2)$ arcs, which contradicts our assumption. 

(iii) As $u\notin V$ and $v^-\in V$, there is an arc between $u$ and $v^-$. Suppose that $u$ dominates $v^-$, then by (i) we have, $v^{--},u^+,v^-\in V$, contradicting that by (ii) we must have $v^{--}v^-\in A(D)$. In the same way we conclude from (i) and (ii) that $v\dom u^+$.
\qed

\begin{lemma}\label{nonrealarc}
	Suppose that (\ref{C*}) holds. If $C^g_{i}, i\in[2]$ has a pair of consecutive vertices $u,u^+$ which  belong to the same partite set $V$, then $C^g_{3-i}$ contains no vertex from $V$. 
\end{lemma}
\pf Suppose the lemma is false and let $w$ be a vertex of $C^g_{3-i}$ which belongs to the partite set $V$. Consider first the case when all vertices of $C^g_{3-i}$ belong to $V$. Then all vertices of  $C^g_{3-i}$ can be inserted between $u$ and $u^+$, which contradicts the assumption (\ref{C*}).
	Hence we may assume that $\ell{}(C^g_{3-i})>0$. 

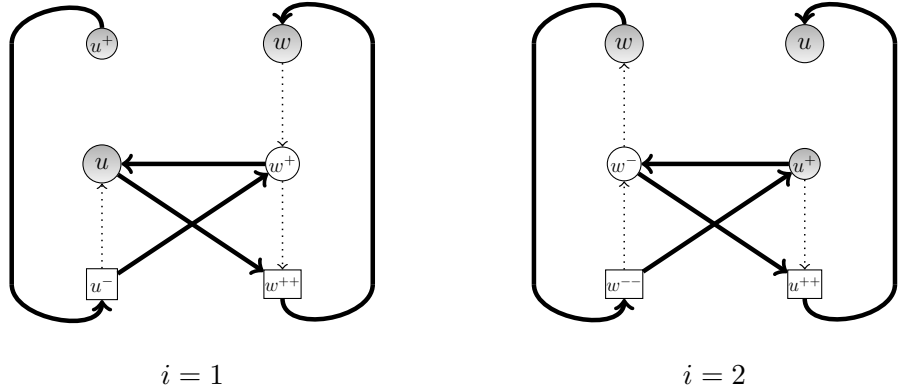
\begin{figure}[H]
\begin{center}
\begin{tabular}{ccc}
\tikzstyle{vertexL}=[circle,draw, minimum size=18pt, scale=0.8, inner sep=0.5pt]
\tikzstyle{vertexS}=[circle,draw, minimum size=18pt, scale=0.65, inner sep=0.5pt]
\tikzstyle{vertexSS}=[circle,draw, minimum size=18pt, scale=0.6, inner sep=0.5pt]
\tikzstyle{vertexLx}=[circle,draw, top color=black!30, bottom color=black!5, minimum size=18pt, scale=0.8, inner sep=0.5pt]
\tikzstyle{vertexSx}=[circle,draw, top color=black!30, bottom color=black!5, minimum size=18pt, scale=0.65, inner sep=0.5pt]
\tikzstyle{vertexSSx}=[circle,draw, top color=black!30, bottom color=black!5, minimum size=18pt, scale=0.5, inner sep=0.5pt]
\tikzstyle{vertexLy}=[rectangle,draw, minimum size=18pt, scale=0.8, inner sep=0.5pt]
\tikzstyle{vertexSy}=[rectangle,draw, minimum size=18pt, scale=0.65, inner sep=0.5pt]
\tikzstyle{vertexSSy}=[rectangle,draw, minimum size=18pt, scale=0.6, inner sep=0.5pt]
\begin{tikzpicture}[scale=0.4]
   \node at (8,-1) {{\small $i=1$}};  
 
  \node (w) at (11,10) [vertexLx]{$w$};
  \node (wp) at (11,6) [vertexS]{$w^+$};
  \node (wpp) at (11,2) [vertexSSy]{$w^{++}$};

  \node (up) at (5,10) [vertexSx]{$u^+$};
  \node (u) at (5,6) [vertexLx]{$u$};
  \node (um) at (5,2) [vertexSy]{$u^-$};

  \draw[->, dotted, line width=0.02cm] (w) to (wp);
  \draw[->, dotted, line width=0.02cm] (wp) to (wpp);

  \draw[line width=0.06cm] (up) to [out=90, in=90]   (2,10);
  \draw[line width=0.06cm] (2,10) to (2,2);
  \draw[->, line width=0.06cm] (2,2) to [out=270, in=270]   (um);

  \draw[line width=0.06cm] (wpp) to [out=270, in=270]   (14,2);
  \draw[line width=0.06cm] (14,2) to (14,10);
  \draw[->, line width=0.06cm] (14,10) to [out=90, in=90]   (w);

  \draw[->, dotted, line width=0.02cm] (um) to (u);

  \draw[->, line width=0.06cm] (um) to (wp);
  \draw[->, line width=0.06cm] (wp) to (u);
  \draw[->, line width=0.06cm] (u) to (wpp);

\end{tikzpicture} & & 
\tikzstyle{vertexL}=[circle,draw, minimum size=18pt, scale=0.8, inner sep=0.5pt]
\tikzstyle{vertexS}=[circle,draw, minimum size=18pt, scale=0.65, inner sep=0.5pt]
\tikzstyle{vertexSS}=[circle,draw, minimum size=18pt, scale=0.6, inner sep=0.5pt]
\tikzstyle{vertexLx}=[circle,draw, top color=black!30, bottom color=black!5, minimum size=18pt, scale=0.8, inner sep=0.5pt]
\tikzstyle{vertexSx}=[circle,draw, top color=black!30, bottom color=black!5, minimum size=18pt, scale=0.65, inner sep=0.5pt]
\tikzstyle{vertexSSx}=[circle,draw, top color=black!30, bottom color=black!5, minimum size=18pt, scale=0.5, inner sep=0.5pt]
\tikzstyle{vertexLy}=[rectangle,draw, minimum size=18pt, scale=0.8, inner sep=0.5pt]
\tikzstyle{vertexSy}=[rectangle,draw, minimum size=18pt, scale=0.65, inner sep=0.5pt]
\tikzstyle{vertexSSy}=[rectangle,draw, minimum size=18pt, scale=0.6, inner sep=0.5pt]
\begin{tikzpicture}[scale=0.4]
   \node at (8,-1) {{\small $i=2$}};   

  \node (w) at (5,10) [vertexLx]{$w$};
  \node (wm) at (5,6) [vertexS]{$w^-$};
  \node (wmm) at (5,2) [vertexSSy]{$w^{--}$};

  \node (u) at (11,10) [vertexLx]{$u$};
  \node (up) at (11,6) [vertexSx]{$u^+$};
  \node (upp) at (11,2) [vertexSSy]{$u^{++}$};

  \draw[->, dotted, line width=0.02cm] (wmm) to (wm);
  \draw[->, dotted, line width=0.02cm] (wm) to (w);
  \draw[line width=0.06cm] (w) to [out=90, in=90]   (2,10);
  \draw[line width=0.06cm] (2,10) to (2,2);
  \draw[->, line width=0.06cm] (2,2) to [out=270, in=270]   (wmm);

  \draw[line width=0.06cm] (upp) to [out=270, in=270]   (14,2);
  \draw[line width=0.06cm] (14,2) to (14,10);
  \draw[->, line width=0.06cm] (14,10) to [out=90, in=90]   (u);

  \draw[->, dotted, line width=0.02cm] (up) to (upp);

  \draw[->, line width=0.06cm] (wmm) to (up);
  \draw[->, line width=0.06cm] (up) to (wm);
  \draw[->, line width=0.06cm] (wm) to (upp);

\end{tikzpicture} \\
\end{tabular}
\end{center}
\caption{The G-cycles created in Lemma~\ref{nonrealarc}. The dark vertices lie in the partite set $V$ and the rectangular vertices lie in the partite set $V'$.
The dotted arcs are not in the created G-cycle.}\label{Lem79}
\end{figure}

Assume that $i=1$. 
As $\ell{}(C^g_2)>0$ we can assume that $w$ is chosen so that $w\in V$ and $w^+\not\in V$. 
In particular $ww^+$ is an arc of $D$ and $u$ and $w^+$ are adjacent. If $u\dom w^+$, then the G-cycle $C^g_1[u^+,u]C^g_2[w^+,w]u^+$ 
has vertex set $V(C^g_1)\cup V(C^g_2)$ and $\ell{}(C^g_1)+\ell{}(C^g_2)$ arcs, contradicting the assumption.  So we must have $w^+\dom u$. 
By Lemma \ref{relationship}, there is a partite set $V'\neq V$ so that $u^-,w^{++}\in V'$ and $D$ contains the arcs $ u^- w^+, w^+u,uw^{++}$ and now 
$C^g_1[u^+,u^-]w^+uC^g_2[w^{++},w]u^+$ is a G-cycle with  vertex set $V(C^g_1)\cup V(C^g_2)$ and at least $\ell{}(C^g_1)+\ell{}(C^g_2)$ arcs
(see Figure~\ref{Lem79}), contradicting the assumption.
This completes the case when $i=1$.

Now assume that $i=2$. 
As $\ell{}(C^g_1)>0$ we can assume that $w$ is chosen so that $w \in V$ and $w^-\not\in V$. 
In particular $w^-w$ is an arc of $D$ and $u^+$ and $w^-$ are adjacent. If $w^- \dom u^+$, then the G-cycle $C^g_1[w,w^-]C^g_2[u^+,u]w$      
has vertex set $V(C^g_1)\cup V(C^g_2)$ and $\ell{}(C^g_1)+\ell{}(C^g_2)$ arcs, contradicting the assumption.  So we must have $u^+ \dom w^-$.
By Lemma \ref{relationship}, there is a partite set $V'\neq V$ so that $w^{--},u^{++}\in V'$ and $D$ contains the arcs $w^{--}u^+,u^+w^-,w^-u^{++}$ and now 
$C^g_{1}[w,w^{--}]u^+ w^- C^g_{2}[u^{++},u]w$ is a G-cycle with vertex set $V(C^g_{1})\cup V(C^g_{2})$ and at least $\ell{}(C^g_{1})+\ell{}(C^g_{2})$ arcs (see Figure~\ref{Lem79}),  contradicting the assumption.
\qed\\

The following corollary can be easily obtained by Lemma  \ref{relationship} (iii).
\begin{corollary}\label{non-adjacent}
	Suppose that (\ref{C*}) holds.  Let $u\in C^g_2$, $v\in C^g_1$ be two vertices in different partite sets, then either $C^g_1[v^+,v]C^g_2[u,u^-]$  or  $C^g_1[v^+,v^-]uvC^g_2[u^+,u^-]$ is a spaning $(v^+,u^-)$-G-path in $D[V(C^g_1\cup C^g_2)]$. Furthermore, the relevant path has at least $\ell{}(C^g_1)+\ell{}(C^g_2)-1$ arcs.
\end{corollary}

The following very important observation on irreducible cycle factors of \smd{}s was shown in \cite{yeoJGT24}.
\begin{theorem}\cite{yeoJGT24}\label{irreducibleCF}
	Let $D$ be a \smd{} and let $F$ be an irreducible  cycle factor. Then we can label the cycles in $F$ by $C_1,\ldots{},C_k$ such that $C_i\simeq> C_j$ for each $i<j$.
\end{theorem}
We now prove that this result extends to irreducible G-cycle factors.

\begin{lemma}\label{irreducibleGCF}
	Let $D$ be a \smd{} and let $F$ be an irreducible  G-cycle factor. Then we can label the cycles in $F$ by $C^g_1,\ldots{},C^g_k$ such that $C^g_i\simeq> C^g_j$ for each $i<j$.
\end{lemma}
\pf By Corollary \ref{2cycleorder}  we have $C^g_i\simeq> C^g_j$ or $C^g_j\simeq> C^g_i$ for every pair of distinct G-cycles.
Suppose there is no ordering of the cycles as claimed. Let $T$ be the tournament on $k$ vertices $v_1,\ldots{},v_k$ and with an arc from $v_i$ to $v_j$ if $C^g_i\simeq> C^g_j$. 
Now $T$ has a cycle and therefore also a 3-cycle $v_iv_jv_kv_i$ (take a shortest cycle). Back in $D$ this means that   $C^g_i\simeq> C^g_j$,
$C^g_j\simeq> C^g_k$ and $C^g_k\simeq> C^g_i$. If all three G-cycles are real cycles (contain no pair of consecutive vertices from the same partite set), then it follows from Theorem \ref{irreducibleCF}, that $F$ does not have the minimum number of cycles. Hence one of the cycles, say, $C^g_j$ contains a pair of consecutive vertices $z,z^+$ from the same partite set $V$. By Lemma \ref{nonrealarc}, neither $C^g_i$ nor $C^g_k$ contains a vertex from $V$ and hence $z,z^+$ are adjacent to every vertex in $V(C^g_i)\cup V(C^g_k)$. By the assumption that $C^g_i\simeq> C^g_j\simeq> C^g_k$, we conclude that $z^+$ has an in-neighbour $u$ on $C^g_i$ and $z$ has an out-neighbour $w$ on $C^g_k$.
If $w^-$ and $u^+$ belong to the same partite set or $w^-\dom u^+$, then $C^g_i[u^+,u]C^g_j[z^+,z]C^g_k[w,w^-]u^+$ is a G-cycle on $V(C^g_i)\cup V(C^g_j)\cup V(C^g_k)$ of length  at least $\ell{}(C^g_i)+\ell{}(C^g_j)+\ell{}(C^g_k)$ (here we use that there is no arc between $z$ and $z^+$ so we loose at most two arcs and gain two or three arcs).
Hence we can assume that $u^+\dom w^-$. By Lemma \ref{relationship} and the fact that $C^g_k\simeq> C^g_i$, this implies that $w^-\dom u^{++}$ and $w^{--}\dom u^+$. Now $C^g_i[u^{++},u]C^g_j[z^+,z]C^g_k[w,w^{--}]u^+w^-u^{++}$ is a G-cycle on $V(C^g_i)\cup V(C^g_j)\cup V(C^g_k)$ of length  at least $\ell{}(C^g_i)+\ell{}(C^g_j)+\ell{}(C^g_k)+1$, contradiction. Hence we have shown that there can be no cycle in $T$, meaning that $T$ is the transitive tournament on $k$ vertices whose vertices can be (re)labelled $v_1,\ldots{},v_k$ such that $v_i\dom v_j$ precisely when $i<j$. Translating this back to $D$ we obtain the desired ordering of the G-cycles.
\qed

\begin{lemma}\label{lem-LRnocommon}
	Let $F=\cup_{i\in[k]}C^g_i$ be an irreducible  G-cycle factor and let $(C^g_1,\ldots{},C^g_k)$ be an ordering of its G-cycles such that $C^g_i\simeq> C^g_j$ for each $i<j$. If $C^g_r$ has a pair of consecutive vertices from the same partite set, then there is no partite set $V$ such that both $\cup_{j=1}^{r-1}V(C^g_j)$ and $\cup_{j=r+1}^kV(C^g_j)$ contain vertices of $V$. 
\end{lemma}
\pf Suppose w.l.o.g that $C_r^g$ has a pair of consecutive vertices  $u,u^+$ from the same partite set and that $w\in C^g_p,z\in C^g_q$ with $p<r<q$ are two vertices belonging to the same partite set $V$. By Lemma \ref{nonrealarc}, all vertices in the partite set containing $u,u^+$ are on the cycle $C^g_r$. Hence $w^-$ and $u^+$ (resp., $u$ and $z^+$) belong to different partite sets. Now there are 4 cases according to the direction of the arcs between $w^-$ and $u^+$, respectively between $u$ and $z^+$.
\begin{itemize}
	\item If both $w^-u^+$ and $uz^+$ are arcs, then $C^g_p[w,w^-]C^g_r[u^+,u]C^g_q[z^+,z]w$ is a G-cycle in $D$ with vertex set $V(C^g_p)\cup V(C^g_r)\cup V(C^g_q)$ and precisely  $\ell{}(C^g_p)+\ell{}(C^g_q)+\ell{}(C^g_r)$ arcs as we loose the two arcs $w^-w,zz^+$ and gain the two arcs $w^-u^+,uz^+$. Note that it follows from Lemma \ref{nonrealarc} that $w^-w,zz^+$ are in fact arcs of $D$.
\end{itemize}

It should be noted that in the following three cases we have $|V(C_r^g)|\geq 3$ since otherwise $u$ and $w^{--}$ are in the same partite set or $u^+$ and $z^{++}$ are in the same partite set by Lemma \ref{relationship} (i), which contradicts the fact that all vertices in the partite set containing $u,u^+$ are on the cycle $C^g_r$.

\begin{itemize}
	\item If $u^+w^-$ and $uz^+$ are arcs, then by Lemma \ref{relationship}, $w^{--}u^+$ and $w^-u^{++}$ are both arcs and hence
	$C^g_p[w,w^{--}]u^+w^-u^{++}C^g_r[u^{++},u]C^g_q[z^+,z]w$ is a G-cycle in $D$ with vertex set $V(C^g_p)\cup V(C^g_r)\cup V(C^g_q)$ and precisely  $\ell{}(C^g_p)+\ell{}(C^g_q)+\ell{}(C^g_r)$ arcs as we loose 4 arcs and gain 4 new ones.
	\item If $w^-u^+$ and $z^+u$ are arcs, then by Lemma \ref{relationship}, $uz^{++}, u^-z^+$ are both arcs and hence
	$C^g_p[w,w^-]C^g_r[u^+,u^-]z^+uz^{++}C^g_q[z^{++},z]w$ is a G-cycle in $D$ with vertex set $V(C^g_p)\cup V(C^g_r)\cup V(C^g_q)$ and precisely  $\ell{}(C^g_p)+\ell{}(C^g_q)+\ell{}(C^g_r)$ arcs as we loose 4 arcs and gain 4 new ones.
	\item If $u^+w^-$ and $z^+u$ are arcs, then using Lemma \ref{relationship} twice we see that $w^{--}u^+,w^-u^{++}$, $uz^{++}, u^-z^+$ are all arcs and then $C^g_p[w,w^{--}]u^+w^-C^g_r[u^{++},u^-]z^+uz^{++}C^g_q[z^{++},z]w$ is a G-cycle in $D$ with vertex set $V(C^g_p)\cup V(C^g_r)\cup V(C^g_q)$ and precisely  $\ell{}(C^g_p)+\ell{}(C^g_q)+\ell{}(C^g_r)$ arcs as we loose 6 arcs and gain 6 new ones.
\end{itemize}

In all four cases we obtain a new G-cycle factor with fewer cycles, contradicting the assumption. \qed


\begin{lemma}\label{loseonearc}
	Let $F=\cup_{i\in[k]}C^g_i$ be an irreducible  G-cycle factor and let $(C^g_1,\ldots{},C^g_k)$ be an ordering of the  G-cycles of $F$ such that $C^g_i\simeq> C^g_j$ for each $i<j$. Suppose that $C^g_p$ and $C^g_q$ with $p<q$ contain vertices of the same partite set, then one of the following statements holds.
	
	(i) There is a spanning G-cycle in $D[V(\cup_{i=p}^q C^g_i)]$ with at least $\sum_{i=p}^q \ell{}(C^g_i)-1$ arcs.
	
	(ii) Each cycle $C_i^g$ with $p\leq i\leq q$ is a real cycle. Further, there are two partite sets of $D$ which cover all vertices of such cycles. 
\end{lemma}


\pf Suppose that $v_p\in C^g_p$ and $v_q\in C^g_q$ are two vertices in the same partite set $V_1$. By Lemma \ref{nonrealarc}, there is no pair of  vertices from $V_1$ which are consecutive on one of the cycles of $F$. In particular, $v_p^-,v_q^+\notin V_1$. By Lemma \ref{nonrealarc} again, for each $i\in[q-p-1]$, there is a vertex $v_{p+i}$ on $C^g_{p+i}$ such that $v_{p+i}$ and $v_{p+i-1}^-$ belong to different partite sets. 


If $v_{q-1}^-$ and $v_q^+$ belong to different partite sets, then one can obtain a $(v_p,v_q)$-G-path satisfying (i) by applying Corollary \ref{non-adjacent} to the $q-p$ pairs of consecutive cycles. It should be noted that for the case that $|V(C_r^g)|=2$ with $p<r<q$ and $v_rv_{r-1}^-, v_{r+1}v_r^-\in A(D)$, then $C_{r-1}^g[v_{r-1},v_{r-1}^-]v_r^-v_rC_{r+1}^g[v_{r+1},v_{r+1}^-]$ is a spanning $(v_{r-1},v_{r+1}^-)$-G-path in $D[V(C_{r-1}^g\cup C_{r}^g\cup C_{r+1}^g)]$. Recall that $v_q^+\notin V_1$. So it suffices to consider the case that $v_{q-1}^-,v_q^+$ belong to the same partite set which is distinct from  $V_1$, say w.l.o.g. this is the partite set $V_2$. 

\begin{figure}[H]
\begin{center}

\tikzstyle{vertexLvw}=[regular polygon,regular polygon sides=4,draw, top color=gray!40, bottom color=gray!20, minimum size=18pt, scale=0.75, inner sep=0.8pt]
\tikzstyle{vertexSvw}=[regular polygon,regular polygon sides=4,draw, top color=gray!45, bottom color=gray!20, minimum size=18pt, scale=0.6, inner sep=0.8pt]
\tikzstyle{vertexLvv}=[circle,draw, top color=gray!10, bottom color=gray!50, minimum size=20pt, scale=0.75, inner sep=0.5pt]
\tikzstyle{vertexL}=[circle,draw, minimum size=20pt, scale=0.75, inner sep=0.5pt]
\tikzstyle{vertexS}=[circle,draw, minimum size=20pt, scale=0.6, inner sep=0.5pt]
\tikzstyle{vertexSx}=[circle,draw, top color=black!30, bottom color=black!5, minimum size=20pt, scale=0.6, inner sep=0.5pt]
\begin{tikzpicture}[scale=0.4]
  \node (vp) at (4,10) [vertexLvv]{$v_p$};
  \node (vpm) at (4,6) [vertexL]{$v_p^-$};
  \node (vpmm) at (4,2) [vertexS]{$v_p^{--}$};

  \node at (3,13) {$C_p$};
  \draw[->, line width=0.06cm] (vpmm) to (vpm);
  \draw[->, line width=0.06cm] (vpm) to (vp);
  \draw[line width=0.06cm] (vp)  to [out=90, in=90]  (1,10);
  \draw[line width=0.06cm] (1,10)  to  (1,2);
  \draw[->, line width=0.06cm] (1,2)  to [out=270, in=270]  (vpmm);

  \node (vpIm) at (8,10) [vertexS]{$v_{p+1}^-$};
  \node (vpI) at (8,6) [vertexS]{$v_{p+1}$};
  \node (vpIp) at (8,2) [vertexS]{$v_{p+1}^{+}$};

  \node at (9,13) {$C_{p+1}$};
  \draw[->, line width=0.06cm] (vpIm) to (vpI);
  \draw[->, line width=0.06cm] (vpI) to (vpIp);
  \draw[line width=0.06cm] (vpIp)  to [out=270, in=270]  (11,2);
  \draw[line width=0.06cm] (11,2)  to  (11,10);
  \draw[->, line width=0.06cm] (11,10)  to [out=90, in=90]  (vpIm);

  \node at (14,6) {$\cdots$};
  \node (vqIIm) at (17,10) [vertexS]{$v_{q-2}^-$};
  \node (vqII) at (17,6) [vertexS]{$v_{q-2}$};
  \node (vqIIp) at (17,2) [vertexS]{$v_{q-2}^{+}$};

  \node at (18,13) {$C_{q-2}$};
  \draw[->, line width=0.06cm] (vqIIm) to (vqII);
  \draw[->, line width=0.06cm] (vqII) to (vqIIp);
  \draw[line width=0.06cm] (vqIIp)  to [out=270, in=270]  (20,2);
  \draw[line width=0.06cm] (20,2)  to  (20,10);
  \draw[->, line width=0.06cm] (20,10)  to [out=90, in=90]  (vqIIm);

  \node (vqIm) at (24,10) [vertexSvw]{$v_{q-1}^-$};
  \node (vqI) at (24,6) [vertexS]{$v_{q-1}$};
  \node (vqIp) at (24,2) [vertexS]{$v_{q-1}^{+}$};

  \node at (25,13) {$C_{q-1}$};
  \draw[->, line width=0.06cm] (vqIm) to (vqI);
  \draw[->, line width=0.06cm] (vqI) to (vqIp);
  \draw[line width=0.06cm] (vqIp)  to [out=270, in=270]  (27,2);
  \draw[line width=0.06cm] (27,2)  to  (27,10);
  \draw[->, line width=0.06cm] (27,10)  to [out=90, in=90]  (vqIm);

  \node (vq) at (31,10) [vertexLvv]{$v_{q}$};
  \node (vqp) at (31,6) [vertexLvw]{$v_{q}^{+}$};

  \node at (32,13) {$C_{q}$};
  \draw[->, line width=0.06cm] (vq) to (vqp);

  \draw[line width=0.06cm] (vqp)  to [out=270, in=270]  (31,2);
  \draw[line width=0.06cm] (31,2)  to [out=270, in=270]  (34,2);
  \draw[line width=0.06cm] (34,2)  to  (34,10);
  \draw[->, line width=0.06cm] (34,10)  to [out=90, in=90]  (vq);
        
  \draw[->, line width=0.06cm] (vpI) to (vpm);
  \draw[->, line width=0.06cm] (vpmm) to (vpI);
  \draw[->, line width=0.06cm] (vpm) to (vpIp);
  \draw[dotted, ->, line width=0.06cm] (vpIm) to (vqII); \draw[line width=0.06cm] (vpIm) to (10.25,9); \draw[->, line width=0.06cm] (14.75,7) to (vqII);
  \draw[->, line width=0.06cm] (vqIIm) to (vqI);
  \draw[dotted, line width=0.06cm] (vqIm) to (vqp);
  \draw[dotted, line width=0.06cm] (vq)  to [out=120, in=60]  (vp);
\end{tikzpicture}
\end{center}
\caption{The situation when $v_p$ and $v_q$ belong to the same partite set $V_1$ and $v^-_{q-1}$ and $v^+_q$ belong to the same partite set $V_2\neq V_1$.}\label{PicX}
\end{figure}
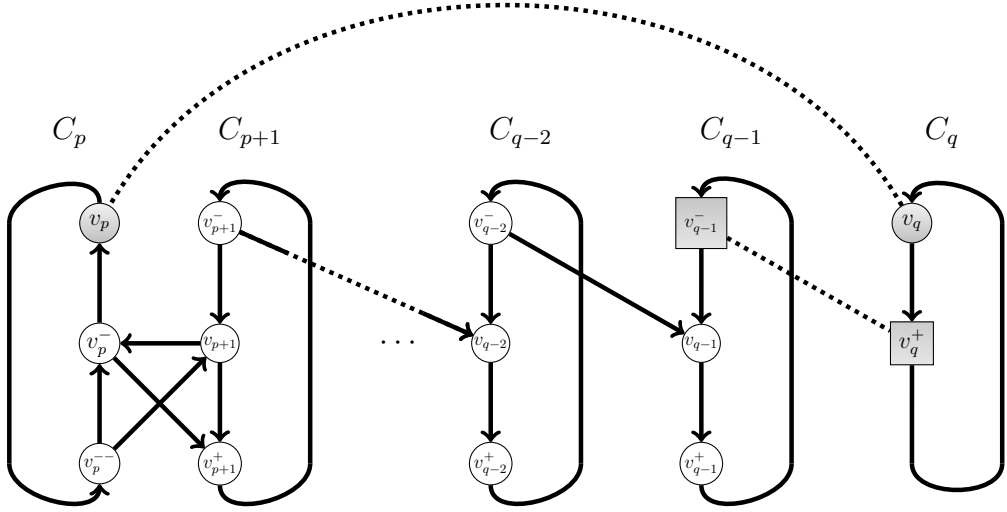

\textbf{Next we show that for each $i\in[q-p-1]$, if $v_{q-i}^-\in V_2$, then $v_{q-i-1}^-\in V_2$ and  $C^g_{q-i}$ is a (real) cycle whose vertices alternate between belonging to  $V_2$ and not belonging to $V_2$.} 

%
%

Recall that $v_{q-1}$ is an arbitrary vertex of $C^g_{q-1}$ which does not belong to the same partite set as $v^-_{q-2}$. For each such vertex we saw that we can assume that $v^-_{q-1}$ is in the same partite set as $v^+_q$. This and Lemma \ref{nonrealarc} implies that every second vertex on $C^g_{q-1}$ belongs to $V_2$. If $v^-_{q-2}\not\in V_2$ then (as we could have chosen the arcs between $v^-_{q-2}$ and $v^-_{q-1}$ instead) we get by the argument above that $v^{--}_{q-1}\in V_2$, contradicting Lemma \ref{nonrealarc}, so we have $v^-_{q-2}\in V_2$. Arguing similarly for the G-cycles $C^g_{q-2},C^g_{q-3},\ldots{},C^g_{p}$ in that order, we conclude that all of the  G-cycles $C^g_{p},\ldots{},C^g_{q-1}$ are cycles of $D$ whose vertices alternate between belonging to $V_2$ and not belonging to $V_2$ and that the vertices  $v^-_p,v^-_{p+1},\ldots{},v^-_{q-1},v^+_q$ all belong to $V_2$.

Now repeating the  derivation above by replacing $(v_p,v_q)$ with $(v_p^-,v_q^+)$  
we conclude that all vertices not in $V_2$ on $C^g_{p+1},\ldots{},C^g_{q-1}$ belong to the same partite set $V_j$ with $j\neq 2$ and $v_p^{--},v_q^{++}$ also belong to $V_j$. That is, the vertices on those cycles alternate between belonging to $V_j$ and belonging to $V_2$. Again, repeating the progress with $(v_p^{--},v_q^{++})$, $(v_p^{---},v_q^{+++})$, $\ldots, (v_p^+,v_q^-)$ in turn, we have that $C_p^g,C_q^g$ are two cycles with vertices  belonging to  $V_2$ and belonging to $V_j$, alternatively. Note that $j=1$ as $v_p,v_q$ are in $V_1$, which shows (ii) holds.\qed

\begin{corollary}\label{coro-lose2arcs}
	If Lemma \ref{loseonearc} (ii) holds, then $C_i^g\compdom C_j^g$ for each $p\leq i<j\leq q$. Moreover, we can obtain a spanning G-cycle $C^g$ in $D[V(\cup_{i=p}^q C^g_i)]$ with $\sum_{i=p}^q \ell{}(C^g_i)-2$ arcs.
\end{corollary}
\pf
Suppose to the contrary that $C_i^g\ncompdom C_j^g$ for some $i,j$ with $p\leq i<j\leq q$. Then there is an arc $uv$ from $C_j^g$ to $C_i^g$. Clearly, $u$ and $v$ belong to different partite sets. W.l.o.g say $u\in V_1$ and $v\in V_2$. Then $u^+\in V_2$ and $v^-\in V_1$, which contradicts Lemma \ref{relationship} (i). So the first statement holds and thus $D[V(\cup_{i=p}^q C^g_i)]$ is non-strong. Further, it is easy to check that $D[V(\cup_{i=p}^q C^g_i)]$ has a spanning G-cycle with $\sum_{i=p}^q \ell{}(C^g_i)-2$ arcs.
\qed

\section{Longest G-cycles}

In this section we consider Problem \ref{prob:longestG-cycle}, that is,  the problem of finding a G-cycle which covers the maximum number of arcs in a \smd{}. The following result shows that we can check in polynomial time whether the given \smd{} has a spanning G-cycle which is also a `cycle' in $D$. This provides the answer to problem \ref{prob:longestG-cycle} when $D$ is a Hamiltonian \smd{}.

\subsection{G-cycles with at least $|V|-k$ arcs for fixed $k$}

\begin{theorem}\label{thm:smdhcalg}\cite{bangJGT29}
	The Hamiltonian cycle problem is solvable in polynomial time for semicomplete multipartite digraphs.
\end{theorem}

For general (non-Hamiltonian) \smd{}s Problem \ref{prob:longestG-cycle} seems challenging as indicated by the results below.
 Even checking whether there is a G-cycle with $|V(D)|-1$ arcs in a non-Hamiltonian \smd{} is already non-trivial, even though we can do as we shall now demonstrate.
A G-cycle with $|V(D)|-1$ arcs is equivalent to a Hamiltonian path in $D$ with both end-vertices in the same partite set.
Now for every vertex $x \in D$ create the \smd{} $D_x$ by adding arcs from $x$ to all other vertices in the same partite set as $x$. 
Now $D_x$ has a Hamiltonian cycle if and only if there is a Hamiltonian path in $D$ ending in $x$ and starting in a vertex from the same partite set as $x$ 
(as $D$ itself was non-Hamiltonian). So there is a G-cycle with $|V(D)|-1$ arcs in $D$ if and only if $D_x$ has a Hamiltonian cycle for some $x \in V(D)$.


Generalizing this idea above we can prove the following.

\begin{theorem}\label{thm:Gcyclen-k$}
	For every fixed $k$ there exists a polynomial algorithm for checking whether a given \smd{} $D$ on $n$ vertices has a spanning G-cycle of length at least $n-k$.
\end{theorem}

\pf Let $D=(V,A)$ be an \smd{}. Note that, by defintion of a G-cycle, a spanning G-cycle $C^g$ of length $n-k$ consists of exactly $k$ vertex-disjoint paths $P_1,P_2,\ldots{},P_k$ and the terminal vertex of $P_i$ belongs to the same partite set as the initial vertex of $P_{i+1}$. Hence if we add all arcs from each of the terminal vertices of $P_1,\ldots{},P_k$ to all other vertices in their partite sets, we obtain a \smd{} with a Hamiltonian cycle containing $P_1,\ldots{},P_k$ as subpaths. For a given subset $X\subset V$ we denote by $D_X$ the \smd{} that we obtain from $D$ by adding all arcs from each vertex $x\in X$ to all other vertices in the same partite set as $x$.
By the remark above, it is not difficult to see that $D_X$ has a Hamiltonian cycle for some choice of $X$ with $|X|=k$ if and only if
$D$ has a spanning G-cycle of length at least $n-k$, where $n=|V|$. There are less than $n^k$ different subsets $X$ to try and together with Theorem \ref{thm:smdhcalg} this  implies that we can check for a spanning G-cycle with at least $n-k$ arcs in polynomial time.
\qed

\begin{problem}
	What is the complexity of checking for given input $k$ and a \smd{} $D$ whether $D$ has a spanning G-cycle of length at least $n-k$, where $n=|V(D)|$?
\end{problem}

\begin{problem}
	Is there an algorithm whose running time is $O(f(k)n^c)$ for some constant $c$ and computable function $f$ (that is, an FPT algorithm) for checking the existence of a G-cycle of length at least $n-k$ in a given \smd{}? 
\end{problem}



\subsection{Bounds of the length of a longest G-cycle (factor)}

Recall that in a strong \smd{}, each maximum G-cycle is spanning (Lemma \ref{longestG-cycle}). The following lemma extends Lemma \ref{longestG-cycle}.
\begin{lemma}
	Let $D$ be a strong \smd{} and let $F$ be any G-cycle subdigraph. Then $D$ has a G-cycle factor with at least $|A(F)|$ arcs.
\end{lemma}
\pf Insert as many vertices into the cycles of $F$ as possible  such that the resulting G-cycle subdigraph still has at least $|A(F)|$ arcs. For convenience, we still call the resulting G-cycle subdigraph $F$. If $F$ is spanning, then there is nothing to prove. Let $v$ be a vertex not covered by $F$. If  there is a G-cycle $C^g\in F$ and  a vertex $u\in C^g$ such that $u$ and $v$ belong to the same partite set or $v$ has both an in-neighbour and an out-neighbour on $C^g$, then it follows just as in the proof of  Lemma \ref{longestG-cycle} that we can insert $v$ in $C^g$ as still have at least the same number of arcs in the new cycle. So we can assume there is no such vertex $v$ in $V(D)-V(F)$.

So for every vertex $v$ in $V(D)-V(F)$, $v$ is adjacent to every vertex in $F$ and moreover for 
each G-cycle $C^g\in F$, either $v\compdom C^g$ or
$C^g\compdom v$. 
We may further assume that no two vertices in $V(D)-V(F)$ belong to the same partite set as the  vertices from the same partite set form a trivial G-cycle which we could add to $F$. This implies that $D-V(F)$ is a transitive tournament. Let  $v$ be  the vertex in $V(D)-V(F)$ which has in-degree zero in $D-V(F)$.

Since $D$ is strong, there is an arc from $F$ to $v$. Assume w.l.o.g. that the arc is from $C^g_1\in F$, then $C^g_1\compdom v$ and hence we must have $C^g_1\compdom D-F$ as no vertex in $V(D-F)$ can be inserted into $F$. Using the fact that $D$ is strong again, $C^g_1$ has an in-arc, say w.l.o.g. it is from $C^g_2$. Then we must have $C^g_2\compdom v$, implying that $C^g_2\compdom D-F$. In the same way, $V(C^g_1)\cup V(C^g_2)$ has an in-arc, say w.l.o.g. from $C^g_3$. Again we conclude that $C^g_3\compdom D-F$. Repeating  this process, we obtain that $F\compdom D-F$,  contradicting the fact that $D$ is strong.  \qed

\bigskip

Let $c_f$ be the maximum number of arcs in a G-cycle factor of $D$. This number and a corresponding
G-cycle factor can be found in polynomial time by Lemma \ref{lem:maxnoarcsCF}.  For any two vertices $x,y\in V(D)$, let $N(x,y)$ be the minimum number of pairs of consecutive non-adjacent vertices  (vertices from the same partite set)
in an $(x,y)$-G-path among all possible $(x,y)$-paths and let $N=\max\{N(x,y):x,y\in V(D)\}$. 

Note that $N=0$ if $D$ is strong.

\bigskip

\textbf{Observation 1:} The maximum number of arcs in a spanning G-cycle of $D$ is at most $\min\{n-N,c_f\}$, where $n=|V(D)|$.

\begin{figure}[H]
\begin{center}
\tikzstyle{vertexLa}=[circle,draw, minimum size=12pt, scale=0.85, inner sep=0.5pt]
\tikzstyle{vertexLb}=[circle, draw, top color=gray!30, bottom color=gray!80, draw, minimum size=12pt, scale=0.85, inner sep=0.5pt]
\tikzstyle{vertexLI}=[diamond, draw, top color=gray!80, bottom color=gray!30,  minimum size=12pt, scale=0.85, inner sep=0.5pt]
\tikzstyle{vertexLII}=[regular polygon,regular polygon sides=4,draw, top color=gray!10, bottom color=gray!60,  minimum size=12pt, scale=0.85, inner sep=0.5pt]
\tikzstyle{vertexLIII}=[star, star points=10, draw, top color=gray!70, bottom color=gray!90, minimum size=12pt, scale=0.85, inner sep=0.5pt]
\begin{tikzpicture}[scale=0.4] 


  \node (x1) at (2,8) [vertexLa]{};
  \node (y1) at (2,2) [vertexLI]{};
  \draw[->, line width=0.06cm] (x1)  to [out=290, in=70]  (y1);
  \draw[->, line width=0.06cm] (y1)  to [out=110, in=250]  (x1);

  \node (x2) at (6,8) [vertexLb]{};
  \node (y2) at (6,2) [vertexLI]{};
  \draw[->, line width=0.06cm] (x2)  to [out=290, in=70]  (y2);
  \draw[->, line width=0.06cm] (y2)  to [out=110, in=250]  (x2);

  \node (x3) at (10,8) [vertexLII]{};
  \node (y3) at (10,2) [vertexLI]{};
  \draw[->, line width=0.06cm] (x3)  to [out=290, in=70]  (y3);
  \draw[->, line width=0.06cm] (y3)  to [out=110, in=250]  (x3);

  \node (x4) at (14,8) [vertexLII]{};
  \node (y4) at (14,2) [vertexLa]{};
  \draw[->, line width=0.06cm] (x4)  to [out=290, in=70]  (y4);
  \draw[->, line width=0.06cm] (y4)  to [out=110, in=250]  (x4);

  \node (x5) at (18,8) [vertexLII]{};
  \node (y5) at (18,2) [vertexLb]{};
  \draw[->, line width=0.06cm] (x5)  to [out=290, in=70]  (y5);
  \draw[->, line width=0.06cm] (y5)  to [out=110, in=250]  (x5);

  \node (x6) at (22,8) [vertexLII]{};
  \node (y6) at (22,2) [vertexLIII]{};
  \draw[->, line width=0.06cm] (x6)  to [out=290, in=70]  (y6);
  \draw[->, line width=0.06cm] (y6)  to [out=110, in=250]  (x6);

  \node (x7) at (26,8) [vertexLa]{};
  \node (y7) at (26,2) [vertexLIII]{};
  \draw[->, line width=0.06cm] (x7)  to [out=290, in=70]  (y7);
  \draw[->, line width=0.06cm] (y7)  to [out=110, in=250]  (x7);

  \node (x8) at (30,8) [vertexLb]{};
  \node (y8) at (30,2) [vertexLIII]{};
  \draw[->, line width=0.06cm] (x8)  to [out=290, in=70]  (y8);
  \draw[->, line width=0.06cm] (y8)  to [out=110, in=250]  (x8);

  \draw[->, line width=0.06cm] (x2)  to  (x1);
  \draw[->, line width=0.06cm] (x3)  to  (x2);
  \draw[->, line width=0.06cm] (y4)  to  (y3);
  \draw[->, line width=0.06cm] (y5)  to  (y4);
  \draw[->, line width=0.06cm] (y6)  to  (y5);
  \draw[->, line width=0.06cm] (x7)  to  (x6);
  \draw[->, line width=0.06cm] (x8)  to  (x7);

\end{tikzpicture}
\end{center}
\caption{A strong \smd{} with a cycle factor consisting of $\frac{n}{2}$ 2-cycles (similar looking vertices belong to the same partite set). All arcs not shown are from left to right. Here $N=0$ and $c_f=n=|V(D)|$ but every spanning G-cycle has at most $n-3$ arcs.}
\label{figbadstrongSMD}
\end{figure}
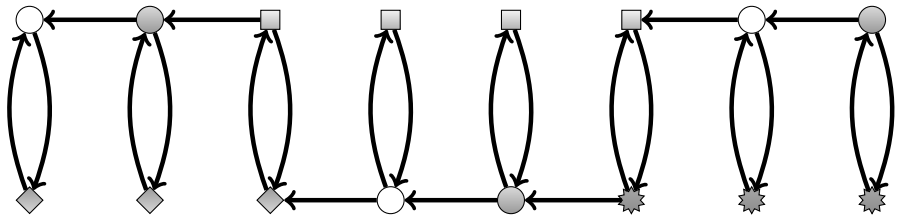

The following examples show that we may not always be able to achieve equality.
\begin{itemize}
	\item Consider a  cycle factor with just two cycles in a non-strong semicomplete bipartite digraph. 
	Then $c_f=n,N=1$ so  $\min\{n-N,c_f\}=n-1$. However, in a bipartite digraph with equally many vertices in each partite set any spanning G-cycle
with consecutive vertices from one partite set must also have consecutive vertices from the other partite set. So there cannot be a 
spanning G-cycle with at least $n-1$ arcs.

	
	\item Consider the  strongly connected \smd{} $D$ in Figure \ref{figbadstrongSMD}, $N=0,c_f=n$, but every spanning G-cycle contains at least three pairs of consecutive vertices which are non-adjacent.
\end{itemize}


Consider the next two results on cycle factors

\begin{theorem} \cite[Theorem 7.9.1]{bang2018}
	If $D$ is a $k$-strong digraph and $X\subseteq V(D)$ has independence number $\alpha(D[X])\leq k$, then there is a cycle subdigraph in $D$ that covers $X$. 
\end{theorem}
\begin{theorem}\cite[Theorem 7.9.2]{bang2018}
	Let $D$ be  a $(\lfloor k/2\rfloor+1)$-strong \smd{} and let  $X$ be an arbitrary set of vertices in $D$ such that $X$ includes at most $k$ vertices from any partite set of $D$. If there exists a cycle subdigraph in $D$ which covers $X$, then there is a cycle $C$ in $D$ such that $X\subseteq V(C)$. 
\end{theorem}

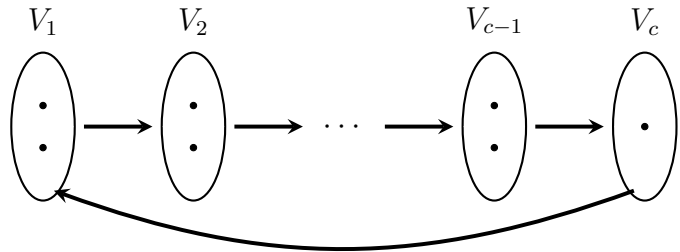
\begin{figure}[H]
	\centering\begin{tikzpicture}[scale=0.4]
	\foreach \i in {(-10,0),(-5,0),(5,0),(10,0)}{\draw[-stealth, line width=0.8pt] \i ellipse [x radius=30pt, y radius=70pt];}
	\coordinate [label=center:$V_1$] () at (-10,3.5);
	\coordinate [label=center:$V_2$] () at (-5,3.5);
	\coordinate [label=center:$\cdots$] () at (0,0);
	\coordinate [label=center:$V_{c-1}$] () at (5,3.5);
	\coordinate [label=center:$V_c$] () at (10,3.5); 
	
	\filldraw[black](-10,0.7) circle (3pt)node[](){};
	\filldraw[black](-10,-0.7) circle (3pt)node[](){};
	\filldraw[black](-5,0.7) circle (3pt)node[](){};
	\filldraw[black](-5,-0.7) circle (3pt)node[](){};
	\filldraw[black](5,0.7) circle (3pt)node[](){};
	\filldraw[black](5,-0.7) circle (3pt)node[](){};
	\filldraw[black](10,0) circle (3pt)node[](){};

	\filldraw[white](-9,0) circle (0.1pt)node[](u){};
	\filldraw[white](-6,0) circle (0.1pt)node[](v){};
	\filldraw[white](-4,0) circle (0.1pt)node[](u1){};
	\filldraw[white](-1,0) circle (0.1pt)node[](v1){};
	\filldraw[white](1,0) circle (0.1pt)node[](u2){};
	\filldraw[white](4,0) circle (0.1pt)node[](v2){};
	\filldraw[white](6,0) circle (0.1pt)node[](u3){};
	\filldraw[white](9,0) circle (0.1pt)node[](v3){};
	\filldraw[white](10,-2) circle (0.1pt)node[](u4){};
	\filldraw[white](-10,-2) circle (0.1pt)node[](v4){};
	\foreach \i/\j in {u/v,u1/v1,u2/v2,u3/v3}{\draw[-stealth, line width=1.5pt] (\i) edge (\j);}
	\draw[-stealth, line width=1.5pt] (u4) edge[bend left=20] (v4);
	\end{tikzpicture}\caption{A $k$-strong \smd{} $D^{\prime}$ with partite sets
		$V_1,\dots{},V_c$ such that $|V_i|=k+1$ for $i\in [c-1]$ and $|V_c|=k$. All arcs between different $V_i$s which are not shown are from left to right.}\label{noclose}
\end{figure}

A corollary of these results is that every $k$-strong \smd{} $D$ with no partite set larger than $k$ has a Hamiltonian cycle \cite{yeoJGT24}. It is natural to ask whether we can lower the connectivity slightly if we just want $D$ to contain a spanning G-cycle with close to $n$ arcs, say $n-f(k)$ arcs for some computable function $f$. This turns out to be wrong. It is easy to check that when $k\geq 1$, every G-cycle factor in the $k$-strong \smd{} $D'$ in Figure \ref{noclose} has at most $n-c+1$ arcs where $c$ is the number of partite sets of $D'$.

Recall from Section \ref{sec:longestGP} that we can find a longest G-path in polynomial time and as the length of a longest G-path is at least the length of a longest G-cycle minus one, we can determine an upper bound on the length of a longest G-cycle in polynomial time by finding a longest G-path. This bound may be far from sharp as indicated by the example in Figure \ref{figbadstrongSMD} as every \smd{} with a cycle factor has a Hamiltonian path by Theorem \ref{thm:longpathSMD} but we can continue the pattern in the figure to make the maximum number of arcs in a G-cycle arbitrarily far from $n$.

%
%
%
%
%
%
%
%
%
%

\subsection{Results in semicomplete bipartite digraphs}

\begin{theorem}\cite{gutinVANB1984,haggkvistC9}
	A semicomplete bipartite digraph has a Hamiltonian cycle if and only if it has a (real) cycle factor and it is strong. 
\end{theorem} 

\begin{theorem}\cite{gutinK1987}
	Let $D$ be a strong semicomplete bipartite digraph. The length of a longest cycle in $D$ is equal to the maximum order of a cycle subdigraph in $D$.
\end{theorem}

\begin{lemma}
	Let $D$ be a semicomplete bipartite digraph with partite sets $V_1$ and $V_2$. Suppose that $F$ is an irreducible G-cycle  factor in $D$. Let $(C^g_1,\ldots{},C^g_k)$  with $k\geq 2$ be an ordering of vertex-disjoint G-cycles in $F$ such that $C^g_i\simeq> C^g_j$ for each $i<j$.  Then $C^g_i\compdom C^g_j$ for each $i<j$. Moreover, either $F=C_1^g\cup C_2^g$ with $V(C_1^g)=V_{t},V(C_2^g)=V_{3-t}$ for some $t\in[2]$ or, each $C_i^g$ is a real cycle and $D[V(F)]$ has a spanning G-cycle with $|A(F)|-2$ arcs.
\end{lemma} 
\pf  If all cycles in $F$ are real cycles, then the vertices of each cycle alternate between belonging to $V_1$ and belonging to $V_2$. In particular, $C_1^g$ and $C_k^g$ contain vertices of the same partite set. It follows by Lemma \ref{loseonearc} and Corollary \ref{coro-lose2arcs} that the claim holds.

So we may assume  that there exists a G-cycle $C_r^g$ in $F$ containing two consecutive vertices $w,w^+$ which belong to the same partite set. W.l.o.g assume that $w,w^+$ belong to $V_1$.  We may further assume that $F$ contains no trivial G-cycle all of whose  vertices belong to $V_1$ since otherwise the vertices of such a trivial G-cycle can be inserted between $w$ and $w^+$, which contradicts the assumption that  $F$ is irreducible.
	
	By Lemma \ref{lem-LRnocommon}, we may assume that $\cup_{i\in[r-1]}V(C_i^g)\subseteq V_1$ and   $\cup_{j\geq r+1}V(C_j^g)\subseteq V_2$. 
Since $F$ contains no trivial G-cycle  all of whose  vertices belong to $V_1$ and $F$ is irreducible, we have that $r=1$, $F=C_1^g\cup C_2^g$ and $V(C_2^g)\subseteq V_2$. Moreover, we have $C^g_1\compdom C^g_2$. If not, assume that $uv$ is an arc from $C^g_2$ to $C^g_1$, then we obtain a contradiction with Lemma \ref{relationship} (ii) as both $u$ and $u^+$ belong to $V_2$. Next we show that $V(C_1^g)\subseteq V_1$. If not, choose a vertex $z\in V_2$ on $C_1^g$ with $z^-\in V_1$. Then all vertices of $C_2^g$ can be inserted between $z^-$ and $z$, which contradicts the assmption that $F$ is irreducible. As $F$ is spanning we conclude that $V(C_1^g)=V_{1}$ and $V(C_2^g)=V_{2}$.\qed

\subsection{Consequences for complete split digraphs}

A {\bf split graph} is a graph whose vertex set can be partitioned into two sets such that one induces  a clique and the other an independent set. Between these two sets there can be arbitrary connections. A {\bf split digraph} is a digraph $D=(V_1,V_2,A)$ where $V(D)$ has a partition $V_1, V_2$ such that the subdigraph induced by $V_1$ has no arcs and the subdigraph induced by $V_2$ is semicomplete, so it is any digraph whose underlying graph is a split graph. A {\bf semicomplete split digraph} is a split digraph in which there is at least one arc between every of $V_1$ and every vertex of $V_2$.

Note that a semicomplete split digraph is a semicomplete multipartite digraph for which at most  one partite set has size more than one. If there is such a set, then without loss of generality this is the set $V_1$.
The following lemma follows from Lemma~\ref{getlongGcycle} with $c^{\prime}=1$.

\begin{lemma}
	Let $D$ be a strong semicomplete split digraph with a G-cycle factor $F$. Then $D$ has a spanning G-cycle with at least $|A(F)|-1$ arcs. 
\end{lemma}

\begin{lemma}\label{getlongGcycle}
	Let $D$ be a strong \smd{} with a G-cycle factor $F$. Then $D$ has a spanning G-cycle with at least $|A(F)|-2c^{\prime}$ arcs, where $c^{\prime}$ is the number of non-trivial partite sets of $D$. When $c^{\prime}=1$ $D$ has a spanning G-cycle with at least $|A(F)|-1$ arcs. Furthermore, we can construct such a cycle in polynomial time from a given G-cycle factor $F$.
	
\end{lemma}
\pf We may w.l.o.g assume that $F$ is irreducible and let $(C^g_1,\ldots{},C^g_k)$ be an ordering of vertex-disjoint G-cycles in $F$ such that
$C^g_i\simeq> C^g_j$ for each $i<j$. Since $D$ is strong, there is a path from $C^g_k$ to $C^g_1$. For each backward arc $uv$, say it is from $C^g_j$ to $C^g_i$ with $i<j$, Lemma \ref{feasible} shows that $u_{C^g_j}^+$ and $v_{C^g_i}^-$ belong to the same partite set. 

Choose $C^g_j$ with maximum index such that $C^g_1$ and $C^g_j$ have vertices in the same partite set. It follows by Lemma \ref{loseonearc} that we can merge $C^g_1,\ldots,C^g_j$ to one G-cycle $\tilde{C}^g_1$ such that we loose at most two arcs, that is,
$\ell{}(\tilde{C}^g_1)\geq \sum_{i\in[j]}\ell{}(C^g_i)-2$. Furthermore, Lemma \ref{loseonearc} shows that if there is only one non-trivial partite set, i.e., $c^{\prime}=1$, then $\ell{}(\tilde{C}^g_1)\geq \sum_{i\in[j]}\ell{}(C^g_i)-1$. It should be noted that $C^g_j$ exists as there is a backward arc entering $C^g_1$. Let $(\tilde{C}^g_1,C^g_{j+1}\ldots{},C^g_k)$ be the new irreducible ordering of vertex-disjoint G-cycles. Choose $C^g_q$ with maximum index $q\geq j+1$ such that $\tilde{C}^g_1$ and $C^g_q$ have vertices in the same partite set and repeat the above process, we can obtain the wanted spanning G-cycle. It is worth noting that the merging can be done in at most $c^{\prime}$ steps as there are only $c^{\prime}$ non-trivial partite sets.
\qed

\begin{corollary}
	If a  strong semicomplete split digraph $D$ on $n$ vertices has a cycle factor, then $D$ has a cycle of length at least $n-1$ and we can find a longest cycle in polynomial time.
\end{corollary} 
\pf Since we can convert any cycle factor to an irreducible cycle factor in polynomial time \cite{yeoJGT24} we can find a cycle of length at least $n-1$ by the algorithm in Lemma \ref{getlongGcycle}. Using Theorem \ref{thm:smdhcalg} we can check in polynomial time whether $D$ also has a Hamiltonian cycle. \qed

For general non-semicomplete split digraphs the complexity of the Hamiltonian cycle problem is open even for split digraphs $D=(V_1,V_2, A)$ where $|V_1|=2$. If $|V_1|=1$, the  Hamiltonian cycle problem is equivalent to deciding the existence of a Hamiltonian path with prescribed starting and ending vertices and hence polynomial (but highly non-trivial) \cite{bangJA13}.

\section{Back to $\{0,1\}$-TSP}

The results on G-cycles in \smd{}s in this paper imply the following for the $\{0,1\}$-weighted version of directed TSP: Let $D=(V,A)$ be a semicomplete digraph with $\{0,1\}$-weights on the arcs so that the arcs of weight 1 induce vertex-disjoint complete subdigraphs $D_1,\ldots{},D_s$ of $D$. Let $D_0$ be the \smd{} that we obtain from $D$ by deleting all arcs of weight 1 from $D$.
\begin{itemize}
\item We can find a minimum cost Hamiltonian path of $D$ in polynomial time. (Theorem \ref{thm:longG-path})
\item If $D_0$ is an extended semicomplete digraph, then we can solve the TSP problem on $D$ in polynomial time. (Theorem \ref{thm:longGCextSD})
  \item For every fixed $k$ there is a polynomial algorithm for checking whether a $\{0,1\}$-arc-weighted semicomplete digraph $D$ where the arcs of cost 1 induce complete digraphs has a TSP tour of weight at most $k$. (Theorem \ref{thm:Gcyclen-k$})
  \end{itemize}

\section{Two NP-complete problems on cycles in semicomplete multipartite digraphs}

 The purpose of this section is to prove two complexity results
   concerning cycles meeting all partite sets in semicomplete multipartite graphs.
 First let us recall that every strong \smd{} has a cycle which visits each partite set at least once. This  follows from Theorem 8.5.6 in \cite{bang2009}. It is also easy to descibe a polynomial algorithm for finding such a cycle. By Theorem \ref{thm:smdhcalg} we can also decide in polynomial time whether $D$ has a Hamiltonian cycle.

The next two results show that, unless P=NP, the following two natural properties of cycles in semicomplete multipartite digraphs  cannot be checked in polynomial time.

\begin{theorem} \label{NPproof1}
	The following problem is NP-complete: Given a semicomplete multipartite digraph $D=(V_1,\ldots{},V_c,A)$ decide whether $D$ contains a cycle $C$ such that $1\leq |V(C)\cap V_i|<|V_i|$ for every $i\in [c]$.
\end{theorem}

\pf  We will reduce from 3-SAT, which is known to be NP-hard.  That is, we ask for a truth assignment so that each clause contains at least one true 
literal. 
The following reduction from 3-SAT to a special path problem has been used in many papers (see e.g. \cite{bangTCS438}): Let $W[u,v,p,q]$ be the digraph (the variable gadget) with vertices $\{u,v,y_0,y_1,y_2,\ldots{},y_p,z_1,z_2,
\ldots,
z_q,z_{q+1}\}$ and the arcs of the two $(u,v)$-paths $u y_0 y_1 \cdots{} y_p v$ and $u z_1 z_2 \cdots{} z_{q+1} v$. 

Let ${\cal F}$ be an instance of 3-SAT with variables $x_1,x_2,\ldots{},x_n$  and clauses $C_1,C_2,\ldots{},C_m$. 
For each variable $x_i$ the ordering of the clauses $C_1,C_2,\ldots{},C_m$ induces an ordering of the occurrences of the literal
$x_i$ 
and the literal $\bar{x_i}$ in these. With each variable $x_i$ we associate a copy of
$W[u_i,v_i,p_i,q_i]$ where the literal $x_i$ occurs $p_i$ times and the literal $\bar{x}_i$ occurs $q_i$ times in the clauses of $\cal F$. 
Let this copy be denoted by $W_i$ and let the vertices in $W_i-\{u_i,v_i\}$ all get a superscript of $i$. 
For each $j\in [p_i]$ the vertex  $y_j^i$  will correspond to the $j$'th  occurence of literal $x_i$ and for each $j\in [q_i]$ the vertex
$z_j^i$ will correspond to the $j$'th  occurence of literal $\bar{x}_i$. Let $V_i^*=\{y_0^i,z_{q_i +1}^i\}$ for all $i \in [n]$.
The vertices in $V_i^*$ do not correspond to any literals.
Identify end-vertices of the above mentioned digraphs by setting $v_i=u_{i+1}$ for $i=1,2,\ldots{},n-1$. Let the resulting digraph be denoted by $D_1$.

\begin{figure}[H]
\begin{center}
\tikzstyle{vertexL}=[circle,draw, top color=gray!30, bottom color=gray!5, minimum size=18pt, scale=0.85, inner sep=0.5pt]
\tikzstyle{vertexLi}=[circle,draw, minimum size=18pt, scale=0.85, inner sep=0.5pt]
\tikzstyle{vertexLii}=[regular polygon,regular polygon sides=4, draw, minimum size=18pt, scale=0.85, inner sep=0.5pt]
\tikzstyle{vertexLiii}=[regular polygon,regular polygon sides=6, draw, minimum size=18pt, scale=0.85, inner sep=0.5pt]

\begin{tikzpicture}[scale=0.37]
  \node (u1) at (1,5) [vertexL]{$u_1$};
  \node (y01) at (4,8) [vertexL]{$y_0^1$};
  \node (y11) at (7,8) [vertexLi]{$y_1^1$};
  \node (y21) at (10,8) [vertexLiii]{$y_2^1$};
  \node (u2) at (13,5) [vertexL]{$u_2$};
  \node (y02) at (17,8) [vertexL]{$y_0^2$};
  \node (y12) at (21,8) [vertexLii]{$y_1^2$};
  \node (u3) at (25,5) [vertexL]{$u_3$};
  \node (y03) at (28,8) [vertexL]{$y_0^3$};
  \node (y13) at (31,8) [vertexLi]{$y_1^3$};
  \node (y23) at (34,8) [vertexLiii]{$y_2^3$};
  \node (v3) at (37,5) [vertexL]{$v_3$};
  
  \node (z11) at (5,2) [vertexLii]{$z_1^1$};
  \node (z21) at (9,2) [vertexL]{$z_2^1$};
  
  \node (z12) at (16,2) [vertexLi]{$z_1^2$};
  \node (z22) at (19,2) [vertexLiii]{$z_2^2$};
  \node (z32) at (22,2) [vertexL]{$z_3^2$};
  
  \node (z13) at (29,2) [vertexLii]{$z_1^3$};
  \node (z23) at (33,2) [vertexL]{$z_2^3$};

  \draw[->, line width=0.03cm] (u1) to (y01);
  \draw[->, line width=0.03cm] (y01) to (y11);
  \draw[->, line width=0.03cm] (y11) to (y21);
  \draw[->, line width=0.03cm] (y21) to (u2);

  \draw[->, line width=0.03cm] (u2) to (y02);
  \draw[->, line width=0.03cm] (y02) to (y12);
  \draw[->, line width=0.03cm] (y12) to (u3);

  \draw[->, line width=0.03cm] (u3) to (y03);
  \draw[->, line width=0.03cm] (y03) to (y13);
  \draw[->, line width=0.03cm] (y13) to (y23);
  \draw[->, line width=0.03cm] (y23) to (v3);

  \draw[->, line width=0.03cm] (u1) to (z11);
  \draw[->, line width=0.03cm] (z11) to (z21);
  \draw[->, line width=0.03cm] (z21) to (u2);

  \draw[->, line width=0.03cm] (u2) to (z12);
  \draw[->, line width=0.03cm] (z12) to (z22);
  \draw[->, line width=0.03cm] (z22) to (z32);
  \draw[->, line width=0.03cm] (z32) to (u3);

  \draw[->, line width=0.03cm] (u3) to (z13);
  \draw[->, line width=0.03cm] (z13) to (z23);
  \draw[->, line width=0.03cm] (z23) to (v3);

\end{tikzpicture}
\end{center}
\caption{An illustration of $D_1$ for the instance ${\cal F}=(x_1 \vee \bar{x}_2 \vee x_3) \wedge
(\bar{x}_1 \vee x_2 \vee \bar{x}_3) \wedge 
(x_1 \vee \bar{x}_2 \vee x_3)$. 
The white circles correspond to clause-1, the white squares to clause-2 and
the white hexagons to clause-3.}\label{Thm92}
\end{figure}
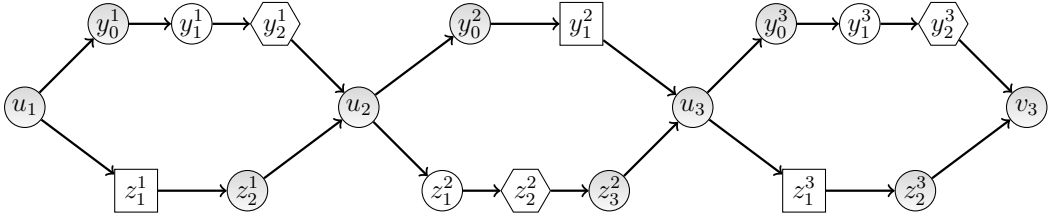

For each $j\in [m]$ we let $V_j$ denote the set of those 3 vertices of $D_1$ which corresponds to the literals of $C_j$ (see Figure~\ref{Thm92}). 
It has been shown several times (see e.g. \cite{bangTCS438}) that  $D_1$ has an $(u_1,v_n)$-path 
which contains at least one vertex from each $V_j$, $j\in [m]$ if and only if $\cal F$ is satisfiable  
(one can analogously show a similar result for a $(u_1,v_n)$-path that avoids at least one vertex from each $V_j$).
Using a similar approach we can in fact show the following.

\2

{\bf Claim A:} $D_1$ has a $(u_1,v_n)$-path which contains exactly one vertex from each $V_i^*$ for all $i\in [n]$ 
and at most two vertices from each $V_j$ for all $j\in [m]$ if and only if $\cal F$ is satisfiable.

\2

{\bf Proof of Claim A:} First assume that $\cal F$ is satisfiable. For every $i \in [n]$ let
$Q_i$ be the path $u_i y_0^i,y_1^i, \ldots, y_{p_i}^i v_i$ if variable $x_i$ is {\em false} and let 
$Q_i$ be the path $u_i z_1^i,z_2^i, \ldots, z_{q_i+1}^i v_i$ otherwise.
The path $Q_1 Q_2 \cdots Q_n$ is now the desired $(u_1,v_1)$-path. 

Conversely, assume that  $D_1$ has a $(u_1,v_n)$-path, $P$, which contains exactly one vertex from each $V_i^*$ for all $i\in [n]$
and at most two vertices from each $V_j$ for all $j\in [m]$.
For every $i \in [n]$, if $u_i y_0^i \in A(P)$ then we let $x_i$ be {\em false} and 
if $u_i z_1^i \in A(P)$ then we let $x_i$ be {\em true}.
As $P$ contains at most two vertices from each $V_j$ ($j\in [m]$) we note that each clause has a true literal (corresponding to the 
vertex not on the path $P$). This completes the proof of Claim~A.\qed

\2

We now extend $D_1$ to a semicomplete multipartite digraph as follows.
Let $V^{start}=\{s_1,s_2\}$ and let $V^{end}=\{t_1,t_2\}$ and add $V^{start} \cup V^{end} \cup \{u^*\}$ to $D_1$ and add all arcs from
$V^{start} \cup \{u^*\}$ to $u_1$ and all arcs from $v_n$ to $V^{end}$. Then add all arcs from  $V^{end}$ to $V^{start} \cup \{u^*\}$.
Let $V^u = \{u^*,u_1,u_2, \ldots, u_{n}, v_n{} \}$. The partite sets in our final digraph, $D$, will be the following.

\[
 V^{start}, V_1, V_2, \ldots, V_m, V_1^*, V_2^*, \ldots, V_n^*, V^u, V^{end}
\]

We now add all arcs from $W_k$ to $W_j$, between vertices of different partite sets, for all $1 \leq j < k \leq n$.
We then add all arcs from $V^{end}$ to $W_i$ and from $W_i$ to $V^{start} \cup \{u^*\}$, 
for all $i \in [n]$.
Then, for all $i \in [n]$ we add all arcs 
from $\{y_0^i,y_1^i, \ldots, y_{p_i}^i,v_i\}$ to $\{u_i,z_1^i,z_2^i, \ldots, z_{q_i+1}^i\}$ between vertices of different partite sets.
Finally for both the paths $u_i y_0^i y_1^i  \cdots y_{p_i}^i v_i$ and $u_i z_1^i z_2^i \cdots z_{q_i+1}^i v_i$ we add all backward arcs 
between non-adjacent vertices. Call the resulting digraph $D$. It is easy to check that $D$ is a semicomplete multipartite digraph.

We will show that $D$ has a cycle such that every partite set contains at least one vertex on the cycle and at least one vertex
not on the cycle if and only if ${\cal F}$ is satisfiable.

\2

First assume that ${\cal F}$ is satisfiable. Let $P$ be the $(u_1,v_n)$-path in $D_1$ that exists due to Claim~A. 
Note that $P$ contains every vertex in $V^u$ except $u^*$ and that the path $P_1 = s_1 P t_1$ does not contain all vertices
from any partite set in $D$. We will now extend $P_1$ to a cycle, $C$, such that $C$ contains at least one, but not all, vertices from
every partite set. Note that $P_1$ contains vertices from every partite set except possibly some partite sets in $\{V_1, V_2, \ldots, V_m\}$.
Let $Q$ denote a set containing exactly one vertex from each partite sets in $\{V_1, V_2, \ldots, V_m\}$ which contain no vertex from $P_1$.
Note that $D[Q]$ is a tournament (all vertices belong to different partite sets) and therefore contains a Hamiltonian path $P_2$. Now
$P_1P_2s_1$ is the desired cycle in $D$.

\2

Conversely assume that there exists a cycle, $C$, in $D$,  such that that every partite set contains at least 
one vertex on the cycle and at least one vertex not on the cycle.
We may without loss of generality assume that $s_1 \in V(C)$ (and $s_2 \not\in V(C)$).
The only out-neighbour of $s_1$ is $u_1$ so $s_1 u_1 \in A(C)$ and $u_1 \in V(C)$. 
This also implies that $u^* \not\in V(C)$ as the only outneighbour of $u^*$ is $u_1$ (and $s_1 u_1 \in A(C)$).
We will now show that for every gadget, $W_i$, either $u_i y_0^i,y_1^i, \ldots, y_{p_i}^i v_i$ or
$u_i z_1^i,z_2^i, \ldots, z_{q_i+1}^i v_i$ is a subpath of $C$. We first consider the case when $i=1$. 

If $u_1 z_1^1 \in A(C)$ then the only arc out of $z_1^1$, which does not go to a vertex in $\{s_1,s_2,u^*,u_1\}$, is
to $z_1^1 z_2^1$. Analogously, the only arc out of $z_2^1$, which does not go to a vertex in $\{s_1,s_2,u^*,u_1,z_1^1\}$, is
to $z_2^1 z_3^1$. Again, the only arc out of $z_3^1$, which does not go to a vertex in $\{s_1,s_2,u^*,u_1,z_1^1,z_2^1\}$, is
to $z_3^1 z_4^1$. Continuing this process we get the path $u_1 z_1^1,z_2^1, \ldots, z_{q_1+1}^1 v_1$.

Now assume that  $u_1 y_0^1 \in A(C)$. For the sake of contradiction assume that 
$u_1 y_0^1,y_1^1, \ldots, y_{p_1}^1 v_1$ is not a subpath of $C$ and let $k$ be the largest integer such that 
$u_1 y_0^1,y_1^1, \ldots, y_{k}^1$ is a subpath of $C$. This implies that the out-neighbour of $y_k^1$ on $C$
must belong to $\{z_1^i,z_2^i, \ldots, z_{q_i+1}^i \}$, say $z_l^1$. However, now it is not difficult to se that 
the path $z_l^1,z_{l+1}^1, \ldots, z_{q_1+1}^1 v_1$ must belong to $C$ (as if we use a backward arc on the path
$u_i z_1^i,z_2^i, \ldots, z_{q_i+1}^i v_i$, then we cannot get the desired cycle $C$). However this contradict the fact 
that $C$ only contains one of the vertices in $V_1^*=\{y_0^1,z_{q_1+1}\}$. 
This contradiction implies that $u_1 y_0^1,y_1^1, \ldots, y_{p_1}^1 v_1$ is a subpath of $C$ as desired.

For $i=2,3, \ldots, n$ (in that order) we can analogously show that the above property holds. We just have to note that if we 
use an arc from some $W_i$ to some $W_j$ with $j<i$ then $C$ cannot contain a vertex from $V^{end}$ (as $u_i$ can only appear on the cycle 
at most once).
So we have now shown that  for every gadget, $W_i$, either $u_i y_0^i,y_1^i, \ldots, y_{p_i}^i v_i$ or
$u_i z_1^i,z_2^i, \ldots, z_{q_i+1}^i v_i$ is a subpath of $C$.
This implies that we have a 
 $(u_1,v_n)$-path in $D_1$ which contains exactly one vertex from each $V_i^*$ for all $i\in [n]$
and at at most two vertices from each $V_j$ for all $j\in [m]$. Therefore, by Claim~A, $\cal F$ is satisfiable.~\qed\\

By Proposition \ref{shortGcycle},  every strong  \smd{} $D=(V_1,\ldots{},V_c,A)$ has a
  G-cycle of length exactly $c$ visiting each $V_i$ exactly once and the proof implies that such a cycle can be constructed in polynomial time. The next result implies that when we replace 'G-cycle' by 'cycle' we obtain an NP-complete problem.

\begin{theorem}  \label{NPproof2}
        The following problem is NP-complete: Given a semicomplete multipartite digraph $D=(V_1,\ldots{},V_c,A)$ decide whether $D$ 
contains a cycle $C$ such that $|V(C)\cap V_i|=1$ for every $i\in [c]$.
\end{theorem}

\pf 
As in the proof of Theorem~\ref{NPproof1} we will reduce from 3-SAT, which is known to be NP-hard.  
We initially follow the same lines as in the proof of Theorem~\ref{NPproof1}. To make reading easier we repeat the first part of that proof.
Let $W[u,v,p,q]$ be the digraph (the variable gadget) with vertices 
$\{u,v,y_0,y_1,y_2,\ldots{},y_p,z_1,z_2, \ldots, z_q,z_{q+1}\}$ and the arcs of the two $(u,v)$-paths $u y_0 y_1 \cdots{} y_p v$ and $u z_1 z_2 \cdots{} z_{q+1} v$.

Let ${\cal F}$ be an instance of 3-SAT with variables $x_1,x_2,\ldots{},x_n$  and clauses $C_1,C_2,\ldots{},C_m$.
For each variable $x_i$ the ordering of the clauses $C_1,C_2,\ldots{},C_m$ induces an ordering of the occurrences of the literal
$x_i$ and the literal $\bar{x_i}$ in these. With each variable $x_i$ we associate a copy of
$W[u_i,v_i,p_i,q_i]$ where the literal $x_i$ occurs $p_i$ times and the literal $\bar{x}_i$ occurs $q_i$ times in the clauses of $\cal F$.
Let this copy be denoted by $W_i$ and let the vertices in $W_i-\{u_i,v_i\}$ all get a superscript of $i$.
For each $j\in [p_i]$ the vertex  $y_j^i$  will correspond to the $j$'th  occurence of literal $x_i$ and for each $j\in [q_i]$ the vertex
$z_j^i$ will correspond to the $j$'th  occurence of literal $\bar{x}_i$.  Let $V_i^*=\{y_0^i,z_{q_i +1}^i\}$ for all $i \in [n]$.
The vertices in $V_i^*$ do not correspond to any literals.
Identify end-vertices of the above mentioned digraphs by setting $v_i=u_{i+1}$ for $i=1,2,\ldots{},n-1$. Let the resulting digraph be denoted by $D_1$.

For each $j\in [m]$ we let $V_j$ denote the set of those 3 vertices of $D_1$ which corresponds to the literals of $C_j$ (see Figure~\ref{Thm92}).
The proof of the following claim is identical to that of Claim~A in Theorem~\ref{NPproof1} and is therefore omitted.

\2

{\bf Claim A:} $D_1$ has a $(u_1,v_n)$-path which contains exactly one vertex from each $V_i^*$ for all $i\in [n]$
and at most two vertices from each $V_j$ for all $j\in [m]$ if and only if $\cal F$ is satisfiable.

\2

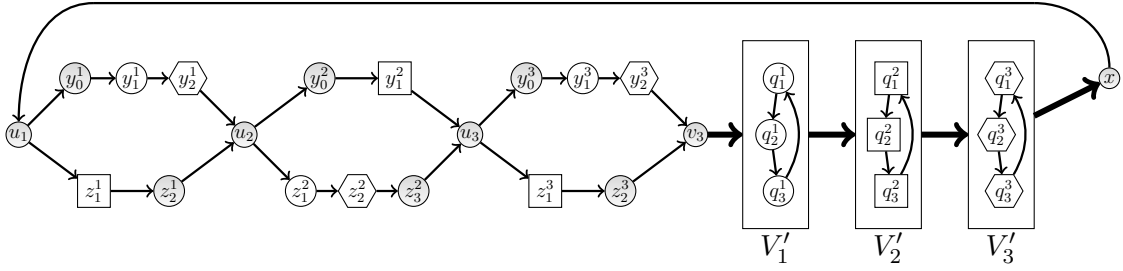
\begin{figure}[H]
\begin{center}
\tikzstyle{vertexL}=[circle,draw, top color=gray!30, bottom color=gray!5, minimum size=12pt, scale=0.65, inner sep=0.5pt]
\tikzstyle{vertexLi}=[circle,draw, minimum size=12pt, scale=0.65, inner sep=0.5pt]
\tikzstyle{vertexLii}=[regular polygon,regular polygon sides=4, draw, minimum size=12pt, scale=0.65, inner sep=0.5pt]
\tikzstyle{vertexLiii}=[regular polygon,regular polygon sides=6, draw, minimum size=12pt, scale=0.65, inner sep=0.5pt]

\begin{tikzpicture}[scale=0.25]
  \node (u1) at (1,5) [vertexL]{$u_1$};
  \node (y01) at (4,8) [vertexL]{$y_0^1$};
  \node (y11) at (7,8) [vertexLi]{$y_1^1$};
  \node (y21) at (10,8) [vertexLiii]{$y_2^1$};
  \node (u2) at (13,5) [vertexL]{$u_2$};
  \node (y02) at (17,8) [vertexL]{$y_0^2$};
  \node (y12) at (21,8) [vertexLii]{$y_1^2$};
  \node (u3) at (25,5) [vertexL]{$u_3$};
  \node (y03) at (28,8) [vertexL]{$y_0^3$};
  \node (y13) at (31,8) [vertexLi]{$y_1^3$};
  \node (y23) at (34,8) [vertexLiii]{$y_2^3$};
  \node (v3) at (37,5) [vertexL]{$v_3$};

  \node (z11) at (5,2) [vertexLii]{$z_1^1$};
  \node (z21) at (9,2) [vertexL]{$z_2^1$};

  \node (z12) at (16,2) [vertexLi]{$z_1^2$};
  \node (z22) at (19,2) [vertexLiii]{$z_2^2$};
  \node (z32) at (22,2) [vertexL]{$z_3^2$};

  \node (z13) at (29,2) [vertexLii]{$z_1^3$};
  \node (z23) at (33,2) [vertexL]{$z_2^3$};

  \node (v1p1) at (41.4,8) [vertexLi]{$q_1^1$};
  \node (v1p2) at (41.0,5) [vertexLi]{$q_2^1$};
  \node (v1p3) at (41.4,2) [vertexLi]{$q_3^1$};
\draw [draw=black] (39.5,0) rectangle (43.0,10);      \node at (41.25,-1) {{\small $V_1'$}};
  \draw[->, line width=0.03cm] (v1p1) to (v1p2);
  \draw[->, line width=0.03cm] (v1p2) to (v1p3);
  \draw[->, line width=0.03cm] (v1p3) to [out=60, in=300]  (v1p1);

  \node (v2p1) at (47.4,8) [vertexLii]{$q_1^2$};
  \node (v2p2) at (47,5) [vertexLii]{$q_2^2$};
  \node (v2p3) at (47.4,2) [vertexLii]{$q_3^2$};
\draw [draw=black] (45.5,0) rectangle (49,10);       \node at (47.25,-1) {{\small $V_2'$}};
  \draw[->, line width=0.03cm] (v2p1) to (v2p2);
  \draw[->, line width=0.03cm] (v2p2) to (v2p3);
  \draw[->, line width=0.03cm] (v2p3) to [out=60, in=300]  (v2p1);

  \node (v3p1) at (53.4,8) [vertexLiii]{$q_1^3$};
  \node (v3p2) at (53,5) [vertexLiii]{$q_2^3$};
  \node (v3p3) at (53.4,2) [vertexLiii]{$q_3^3$};
\draw [draw=black] (51.5,0) rectangle (55,10);    \node at (53.25,-1) {{\small $V_3'$}};
  \draw[->, line width=0.03cm] (v3p1) to (v3p2);
  \draw[->, line width=0.03cm] (v3p2) to (v3p3);
  \draw[->, line width=0.03cm] (v3p3) to [out=60, in=300]  (v3p1);

  \node (x) at (59,8) [vertexL]{$x$};

  \draw[->, line width=0.08cm] (v3) to (39.5,5);
  \draw[->, line width=0.08cm] (43,5) to (45.5,5);
  \draw[->, line width=0.08cm] (49,5) to (51.5,5);
  \draw[->, line width=0.08cm] (55,6) to (x);

  \draw[line width=0.03cm] (x)  to [out=90, in=0]  (55,12);
  \draw[line width=0.03cm] (55,12) to (5,12);
  \draw[->,line width=0.03cm] (5,12)  to [out=180, in=90]  (u1);

  \draw[->, line width=0.03cm] (u1) to (y01);
  \draw[->, line width=0.03cm] (y01) to (y11);
  \draw[->, line width=0.03cm] (y11) to (y21);
  \draw[->, line width=0.03cm] (y21) to (u2);

  \draw[->, line width=0.03cm] (u2) to (y02);
  \draw[->, line width=0.03cm] (y02) to (y12);
  \draw[->, line width=0.03cm] (y12) to (u3);

  \draw[->, line width=0.03cm] (u3) to (y03);
  \draw[->, line width=0.03cm] (y03) to (y13);
  \draw[->, line width=0.03cm] (y13) to (y23);
  \draw[->, line width=0.03cm] (y23) to (v3);

  \draw[->, line width=0.03cm] (u1) to (z11);
  \draw[->, line width=0.03cm] (z11) to (z21);
  \draw[->, line width=0.03cm] (z21) to (u2);

  \draw[->, line width=0.03cm] (u2) to (z12);
  \draw[->, line width=0.03cm] (z12) to (z22);
  \draw[->, line width=0.03cm] (z22) to (z32);
  \draw[->, line width=0.03cm] (z32) to (u3);

  \draw[->, line width=0.03cm] (u3) to (z13);
  \draw[->, line width=0.03cm] (z13) to (z23);
  \draw[->, line width=0.03cm] (z23) to (v3);

\end{tikzpicture}
\end{center}
\caption{An illustration of $D_2$ for the instance ${\cal F}=(x_1 \vee \bar{x}_2 \vee x_3) \wedge 
(\bar{x}_1 \vee x_2 \vee \bar{x}_3) \wedge 
(x_1 \vee \bar{x}_2 \vee x_3)$. The white circles correspond to clause-1, the white squares to clause-2 and
the white hexagons to clause-3. Furthermore, $V_1^*=\{y_0^1,z_2^1\}$, $V_2^*=\{y_0^2,z_{3}^2\}$, $V_2^*=\{y_0^3,z_{2}^3\}$ and
$Q_1^1=\{y_1^1,q_1^1\}$, $Q_2^1=\{z_1^2,q_2^1\}$, $Q_3^1=\{y_1^3,q_3^1\}$,
$Q_1^2=\{z_1^1,q_1^2\}$, $Q_2^2=\{y_1^2,q_2^2\}$, $Q_3^2=\{z_1^3,q_3^2\}$,
$Q_1^3=\{y_2^1,q_1^3\}$, $Q_2^3=\{z_2^2,q_2^3\}$ and $Q_3^3=\{y_2^3,q_3^3\}$.}\label{Thm93}
\end{figure}

We first extend $D_1$ to a new digraph, $D_2$, as follows (this part is different to the proof of Theorem~\ref{NPproof1}).
For all $j \in [m]$ let $V_j'$ be a set of $3$ vertices and let $Q_1^j$, $Q_2^j$ and $Q_3^j$ partition $V_j \cup V_j'$ such that each of
the three sets contain one vertex from $V_j$ and one vertex from $V_j'$ (for example if $V_j=\{a,b,c\}$ and $V_j'=\{d,e,f\}$ then we can let
$Q_1^j=\{a,d\}$, $Q_2^j=\{b,e\}$ and $Q_3^j=\{c,f\}$). Let $D_2$ be the digraph obtained from $D_1$ by first adding all arcs from 
$v_n$ to $V_1'$ and adding all arcs from $V_i'$ to $V_{i+1}'$ for all $i \in [m-1]$. We then add a vertex $x$ and all arcs $V_m'$ to $x$ and
the arc $x u_1$. 
Finally we add a $3$-cycle on the vertices in $V_j'$ for all $j \in [m]$. See Figure~\ref{Thm93} for an example of $D_2$.
Let $U_i=\{u_i\}$ for all $i \in [n]$ and let $U_{n+1}=\{v_n\}$ and let $X=\{x\}$. 
The following claim now holds.

\2

\2

{\bf Claim B:} $D_2$ contains a cycle, $C$, such that the following holds if and only if $\cal F$ is satisfiable.

\begin{itemize}
\item $|U_i \cap V(C)|=1$ for all $i \in [n+1]$.
\item $|V_i^* \cap V(C)|=1$ for all $i \in [n]$.
\item $|Q_i^j \cap V(C)|=1$ for all $i \in [3]$ and $j \in [m]$.
\item $|X \cap V(C)|=1$.
\end{itemize}
 
\2

{\bf Proof of Claim B:} First assume that $\cal F$ is satisfiable. By Claim~A, 
 $D_1$ has a $(u_1,v_n)$-path, $P_1$, which contains exactly one vertex from each $V_i^*$ for all $i\in [n]$
and at most two vertices from each $V_j$ for all $j\in [m]$.
Let $R \subseteq V_1' \cup V_2' \cup \cdots \cup V_m' \cup \{x\}$ be defined such that 
$|(R \cup V(P_1)) \cap Q_i^j|=1$ for all $i \in [3]$ and $j \in [m]$ and $x \in R$.
Note that there is a 
Hamilton path, $P_2$, in $D_2[R]$ starting in a vertex in $V_1'$ and ending in $x$.
Adding the arc $xu_1$ to the path $P_1P_2$ now gives us the desired cycle.

Conversely, assume that  $D_2$ has a cycle $C$ satisfying the conditions of Claim~B.
Note that $\{u_1,v_n,x\} \subseteq V(C)$ (as $|U_1 \cap V(C)|=1$, $|U_{n+1} \cap V(C)|=1$ and $|X \cap V(C)|=1$).
Furthermore the path $C[v_n,x]$ must contain at least one vertex from each $V_j'$, where $j \in [m]$, 
which implies that the path $C[u_1,v_n]$ does not contain all vertices in $V_j$ for any $j \in [m]$.
By Claim~A (and the path $C[u_1,v_n]$) we note that $\cal F$ is satisfiable, which completes the proof of Claim~B.\qed

\2

We now extend $D_2$ to a semicomplete multipartite digraph, $D_3$, as follows.
The partite sets of $D_3$ will be the following.

\begin{itemize}
\item $U_i$ for all $i \in [n+1]$.
\item $V_i^*$ for all $i \in [n]$. 
\item $Q_i^j$ for all $i \in [3]$ and $j \in [m]$.
\item $X$.
\end{itemize} 
 
Note that the above sets partition $V(D_2)$ ($=V(D_3)$).
Firstly, for all $i \in [n]$ we add all arcs
from $\{y_0^i,y_1^i, \ldots, y_{p_i}^i,v_i\}$ 
to $\{u_i,z_1^i,z_2^i, \ldots, z_{q_i+1}^i\}$ between vertices of different partite sets.
Then, for both of the paths $u_i y_0^i y_1^i  \cdots y_{p_i}^i v_i$ and $u_i z_1^i z_2^i \cdots z_{q_i+1}^i v_i$ we add all backward arcs.
We then consider the order $(W_1,W_2,W_3, \ldots, W_n, V_1', V_2', V_3', \ldots, V_m', \{x\})$ and add all arcs 
from sets later in the order to sets earlier in the order, between vertices not in the same partite set. 
That is, if $R$ and $T$ are two sets in the order and $R$ appears before $T$ then we add all arcs from $T$ to $R$, 
between vertices not in the same partite set. 
This completes the construction of $D_3$.

We will show that $D_3$ has a cycle such that every partite set contains exactly one vertex on the cycle if and only if ${\cal F}$ is satisfiable.

\2

First assume that ${\cal F}$ is satisfiable. Let $C$ be the cycle in $D_2$ that exists due to Claim~B.
This cycle is the desired cycle in $D_3$, so we are done in this case.

\2

Conversely assume that there exists a cycle, $C$, in $D_3$, such that that every partite set contains exactly
one vertex on the cycle.
As in the proof of Theorem~\ref{NPproof1} we will show that  for every gadget, $W_i$, either $u_i y_0^i,y_1^i, \ldots, y_{p_i}^i v_i$ or
$u_i z_1^i,z_2^i, \ldots, z_{q_i+1}^i v_i$ is a subpath of $C$. We first consider the case when $i=1$.

If $u_1 z_1^1 \in A(C)$ then the only arc out of $z_1^1$, which does not go to a vertex already visited is $z_1^1 z_2^1$. 
Analogously, the only arc out of $z_2^1$, which does not go to a vertex already visited is $z_2^1 z_3^1$. Continuing this 
process we get the path $u_1 z_1^1,z_2^1, \ldots, z_{q_1+1}^1 v_1$.

Now assume that  $u_1 y_0^1 \in A(C)$. For the sake of contradiction assume that
$u_1 y_0^1,y_1^1, \ldots, y_{p_1}^1 v_1$ is not a subpath of $C$ and let $k$ be the largest integer such that
$u_1 y_0^1,y_1^1, \ldots, y_{k}^1$ is a subpath of $C$. This implies that the out-neighbour of $y_k^1$ on $C$
must belong to $\{z_1^i,z_2^i, \ldots, z_{q_i+1}^i \}$, say it is $z_l^1$. However, now it is not difficult to see that
the path $z_l^1,z_{l+1}^1, \ldots, z_{q_1+1}^1 v_1$ must belong to $C$ (as if we use a backward arc on the path
$u_i z_1^i,z_2^i, \ldots, z_{q_i+1}^i v_i$, then we cannot get the desired cycle $C$). However this contradicts the fact
that $C$ only contains one of the vertices in $V_1^*=\{y_0^1,z_{q_1+1}\}$.
This contradiction implies that $u_1 y_0^1,y_1^1, \ldots, y_{p_1}^1 v_1$ is a subpath of $C$ as desired.

For $i=2,3, \ldots, n$ (in that order) we can analogously show that the above property holds. We just have to note that if we
use an arc from some $W_i$ to some $W_j$ with $j<i$ then $C$ cannot contain the vertex $v_n$.
So we have now shown that  for every gadget, $W_i$, either $u_i y_0^i,y_1^i, \ldots, y_{p_i}^i v_i$ or
$u_i z_1^i,z_2^i, \ldots, z_{q_i+1}^i v_i$ is a subpath of $C$.

The above paths give us the path $C[u_1,v_n]$, which we note is a path in $D_1$. 
The path $C[v_n,u_1]$ must include the vertex $x$ and therefore it must include
at least one vertex from each set $V_j'$ for $j \in [m]$. This implies that the path 
$C[u_1,v_n]$ cannot include all vertices from any set $V_j$, $j \in [m]$. 
By Claim~A we note that ${\cal F}$ is satisfiable, which completes the proof.~\qed

\noindent{}{\bf Conflict of interest statement:}\\
          There are no sources of conflict of interest regarding this paper.
          
          \noindent{}{\bf Data availability statement:}\\
          Data sharing is not applicable to this article as no data sets were generated or analyzed during the
          current study.
          
          \noindent{}{\bf Acknowledgments:}\\
          Financial support from the Independent Research Fund Denmark under grant
          DFF-7014-00037B is gratefully acknowledged. The second author was supported by China Scholarship
          Council (CSC) No. 202106220108.


\end{document}